\documentclass{amsart}
\usepackage{latexsym}
\usepackage{amsmath, amssymb}
\usepackage[all]{xy}
\newtheorem{dfs}{Definition}[section]
\newtheorem{lms}[dfs]{Lemma}
\newtheorem{thms}[dfs]{Theorem}
\newtheorem{props}[dfs]{Proposition}
\newtheorem{qus}[dfs]{Question}
\newtheorem{cors}[dfs]{Corollary}

\newtheorem{conjs}[dfs]{Conjecture}
\newtheorem{exs}[dfs]{Example}

\begin{document}

\title{Flat dimension growth for $C^*$-algebras}
\author{Andrew S. Toms}
\address{Department of Mathematics and Statistics, University of New Brunswick, Fredericton, New Brunswick, Canada, E3B 5A3}
\email{atoms@unb.ca}
\keywords{Nuclear $C^*$-algebras, non-commutative dimension}
\subjclass[2000]{Primary 46L35, Secondary 46L80}
\thanks{This work was supported by an NSERC Discovery Grant}



\begin{abstract}
Simple and nuclear $C^*$-algebras which fail to absorb the Jiang-Su algebra tensorially have 
settled many open questions in the theory of nuclear $C^*$-algebras, but have been 
little studied in their own right.  This is due partly to a dearth of invariants 
sensitive to differences between such algebras.  
We present two new real-valued invariants to fill this void:
the dimension-rank ratio (for unital AH algebras), and the radius of
comparison (for unital and stably finite algebras).
We establish their basic properties, show that they have natural connections to 
ordered $\mathrm{K}$-theory, and prove that the range of the dimension-rank 
ratio is exhausted by simple algebras (this last result shows the class of
simple, nuclear and non-$\mathcal{Z}$-stable $C^*$-algebras to be uncountable).
In passing, we establish a theory of moderate dimension growth for AH algebras,
the existence of which was first supposed by Blackadar.
The minimal instances of both invariants are shown to coincide 
with the condition of being tracially AF among simple unital AH algebras of real 
rank zero and stable rank one, whence they may be thought of as generalised measures of dimension growth.
We argue that the radius of comparison 
may be thought of as an abstract version of the dimension-rank
ratio.
\end{abstract}

\maketitle




\section{Introduction}
The Jiang-Su algebra $\mathcal{Z}$ is by now well known in the
study of nuclear $C^*$-algebras.  All evidence indicates that the
property of being $\mathcal{Z}$-stable --- a $C^*$-algebra $A$ is said to
be $\mathcal{Z}$-stable if $A \otimes \mathcal{Z} \cong A$ ---
is connected naturally to Elliott's program to classify separable and nuclear 
$C^*$-algebras (\cite{El2});  examples due to R{\o}rdam and the author 
show that the largest class of simple, 
separable, unital and nuclear $C^*$-algebras which may be classified up to 
$*$-isomorphism by the Elliott invariant consists of those algebras which are, in addition, 
$\mathcal{Z}$-stable (\cite{R3}, \cite{To1}, \cite{To2}). 
It has been surprising to find that almost all of our stock-in-trade
simple, separable, and nuclear $C^*$-algebras are $\mathcal{Z}$-stable (\cite{JS1}, \cite{TW2}).

Little is known about non-$\mathcal{Z}$-stable $C^*$-algebras in general, save that
they seem able to exhibit arbitrarily strange behaviour.  Specific
examples of such algebras have, over the past several years, been used to settle many open
questions in the theory of separable and nuclear $C^*$-algebras --- see
\cite{R3}, \cite{R5}, \cite{To1}, \cite{To2}, \cite{To3}, \cite{V1}, and \cite{V2} --- but
no attempt has been made to study their structure systematically.  In this paper --- a 
sequel to \cite{To4} in design, though technically independent of it --- we study these algebras
through the introduction of invariants which distill purely non-$\mathcal{Z}$-stable
information:  they are insensitive to differences between $\mathcal{Z}$-stable algebras,
while detecting differences between non-$\mathcal{Z}$-stable $C^*$-algebras 
not readily manifest in known invariants.

In the early sections of the sequel we concentrate on approximately homogeneous (AH) $C^*$-algebras,
as these provide the most tractable examples of simple and non-$\mathcal{Z}$-stable $C^*$-algebras. 
Recall that a {\it homogeneous} $C^*$-algebra has the form
\[
p(\mathrm{C}(X) \otimes \mathcal{K})p,
\]
where $X$ is a compact Hausdorff space, 
$\mathcal{K}$ is the algebra of compact operators on a
separable and infinite-dimensional Hilbert space $\mathcal{H}$, 
and $p \in \mathrm{C}(X) \otimes \mathcal{K}$ is a projection
of constant rank.  
A {\it semi-homogeneous} $C^*$-algebra is a finite direct sum of
homogeneous $C^*$-algebras.

\begin{dfs}[Blackadar, \cite{Bl2}]\label{ahdef} 
An approximately homogeneous (AH) $C^*$-algebra is an inductive limit
\[
A = \lim_{i \to \infty}(A_i,\phi_i),
\]
where each $A_i$ is semi-homogeneous.
\end{dfs}

Let 
\begin{equation}\label{stddecomp}
A \cong \lim_{i \to \infty}(A_i,\phi_i)
\end{equation}
be an unital (i.e., both $A_i$ and $\phi_i:A_i \to A_{i+1}$ are unital for every 
$i \in \mathbb{N}$) AH algebra, where
\begin{equation}\label{stdAi}
A_i := \bigoplus_{l=1}^{m_i} p_{i,l}(\mathrm{C}(X_{i,l}) \otimes \mathcal{K}) p_{i,l}
\end{equation}
for compact Hausdorff spaces $X_{i,l}$, projections 
$p_{i,l} \in \mathrm{C}(X_{i,l}) \otimes \mathcal{K}$, and natural numbers $m_i$.  Put
\[
\phi_{ij} = \phi_{j-1} \circ \phi_{j-2} \circ \cdots \circ \phi_i,
\]
and write $\phi_{i\infty}:A_i \to A$ for the canonical map.
We refer to this collection of objects and maps as a \emph{decomposition} for $A$.
If the $\phi_i$ are injective, then we will refer to this collection as an
{\it injective decomposition}. 

\begin{dfs}\label{flatdim}
Let $A$ be an unital AH algebra.  Say that $A$ has flat dimension growth if
it admits a decomposition for which
\begin{equation}\label{moddimgrowth}
\limsup_{i \to \infty} \ \mathrm{max}_{1 \leq l \leq m_i} \left\{ 
\frac{\mathrm{dim}(X_{i,l})}{\mathrm{rank}(p_{i,l})} \right\} < \infty.
\end{equation}
\end{dfs}

\noindent
A simple unital AH algebra $A$ admitting a decomposition for which (\ref{moddimgrowth})
is zero is said to have {\it slow dimension growth} (\cite{BDR}).  (There are definitions of slow dimension
growth for non-simple algebras in \cite{MP} and \cite{Go1}, but we will not require them here.  Suffice
it to say that these definitions coincide with Definition \ref{flatdim} for simple algebras.)
If (\ref{moddimgrowth}) is finite for some decomposition of $A$, then we may, by passing to a subsequence,
replace the $\limsup$ by a limit;  proofs in the sequel with exploit this.

The beginnings of Definition \ref{flatdim} are contained in Blackadar's 1991 survey article
``Matricial and Ultramatricial Topology'' (\cite{Bl2}).  
At the time, all known simple unital AH algebras had slow dimension growth,
but Blackadar mused nonetheless about the possible existence of a theory of ``AH algebras
with moderate dimension growth'' (synonymous with our flat dimension growth).  His hoped-for theory
was made plausible when Villadsen provided the first examples of simple unital 
AH algebras without slow dimension growth in 1996 (\cite{V1}).
In the sequel we prove that there does indeed exist a theory of flat dimension growth for AH algebras,
and that the natural way to study this theory is through an invariant we call the {\it dimension-rank ratio}.
This invariant for unital AH algebras takes values in the nonnegative reals and recovers, roughly, the minimum possible
value of the limit in (\ref{moddimgrowth}).  It turns out to be of the variety we
seek: it is insensitive to differences between $\mathcal{Z}$-stable algebras  (provided
that they are simple and of real rank zero); it detects subtle differences between
non-$\mathcal{Z}$-stable algebras.  It is also naturally connected to ordered 
$\mathrm{K}$-theory.  These connections lead us
to define an invariant for general unital and stably finite $C^*$-algebras ---
the {\it radius of comparison} --- which measures of the failure of comparison in the
Cuntz semigroup.  This, we argue, is the appropriate abstraction of the
dimension-rank ratio.  Both invariants can be viewed as measuring the ratio
of the matricial size of a $C^*$-algebra to its topological dimension (as constituted
by Kirchberg and Winter's decomposition rank --- see \cite{KW}), despite
the fact that both quantities are frequently infinite for non-$\mathcal{Z}$-stable
algebras.

Our paper is organised as follows:
in section 2 we give the precise definition of the dimension-rank ratio, give a formula
for it in the case of semi-homogeneous $C^*$-algebras, and examine its behaviour with respect to
common constructions;  in section 3 we draw connections between the dimension-rank ratio and
ordered $\mathrm{K}$-theory;  section 4 shows that, among simple algebras of real rank zero,
the minimal instance of the dimension-rank ratio coincides with the condition of being tracially AF;
the range of the dimension-rank ratio is shown to be exhausted by simple algebras
in section 5;  section 6 introduces the radius of comparison, 
defined for any unital and stably finite $C^*$-algebra, and establishes analogues 
of some of our earlier results on the dimension-rank ratio.

\vspace{2mm}
\noindent
{\it Acknowledgements.}
The author would like to thank Bruce Blackadar and Francesc Perera for several helpful discussions, and
Cornel Pasnicu for several comments on earlier versions of this paper.  
Some of this work was done during a workshop at Oberwolfach in August,
2005. The author would like to thank the organisers of the workshop, the institute, 
and its staff for their support.  Finally, we thank the referee for a careful 
reading of our manuscript, and in particular for suggestions which considerably
improved our expsotion in section 2.

\section{The dimension-rank ratio of an AH algebra}

\begin{dfs}\label{drr}
Let $A$ be an unital AH algebra.  Define the dimension-rank ratio of $A$
(write $\mathrm{drr}(A)$) to be the infimum of the set of strictly positive reals $c$ such that
$A$ has a decomposition satisfying
\[
\limsup_{i \to \infty} \mathrm{max}_{1 \leq l \leq m_i}\left\{\frac{\mathrm{dim}(X_{i,l})}{\mathrm{rank}(p_{i,l})}\right\} = c,
\]
whenever this set is not empty, and $\infty$ otherwise.
\end{dfs}

\noindent
By compressing the inductive sequence decomposition for $A$ if necessary, one can replace the
$\limsup$ of Definition \ref{drr} with a limit.
It is sketched in \cite{Bl2} and proved in \cite{Go1} 
that the spaces $X_{i,l}$ in an injective decomposition for
an unital AH algebra $A$ can always be replaced by CW-complexes $\tilde{X}_{i,l}$ of the
same dimension.  From here on we will assume, unless otherwise noted, that the $X_{i,l}$s 
are CW-complexes.  It is also true that if one has a decomposition for $A$ as in Definition
\ref{drr} which is not injective, then it can be replaced with an injective decomposition
for which the limit in Definition \ref{drr} is no larger (cf. \cite{EGL3}).  Thus, we
may assume that the decomposition of the definition is injective whenever this is convenient.

Our first proposition collects some basic properties of the dimension-rank ratio.

\begin{props}\label{drrprop}
Let $A$, $B$ be unital AH algebras, and $I$ an ideal of $A$.  Then:
\begin{enumerate}
\item[(i)] $\mathrm{drr}(A/I) \leq \mathrm{drr}(A)$;
\item[(ii)] $\mathrm{drr}(A \oplus B) = \mathrm{max} \{\mathrm{drr}(A),\mathrm{drr}(B)\}$;
\item[(iii)] $\mathrm{drr}(A \otimes \mathrm{M}_k) \leq (1/k)\mathrm{drr}(A)$;
\item[(iv)] if $A$ and $B$ are simple and of finite dimension-rank ratio, then
$A \otimes B$ has slow dimension growth and $\mathrm{drr}(A \otimes B) = 0$.
\end{enumerate}
\end{props}

\begin{proof} 
For (i), let $\epsilon > 0$ be given, and fix an injective decomposition for $A$ such that for every $i \in \mathbb{N}$
one has
\[
\mathrm{max}_{1 \leq l \leq m_i}\left\{\frac{\mathrm{dim}(X_{i,l})}{\mathrm{rank}(p_{i,l})}\right\} = \mathrm{drr}(A) +\epsilon.
\]
Let $I$ be an ieal of $A$, and let $\phi_{i \infty}:A_i \to A$
be the canonical map.  Then, $I_i := \phi^{-1}_{i \infty}(I)$ is an ideal of $A_i$ for
every $i \in \mathbb{N}$.  Define $\psi_i:A_i/I_i \to A_{i+1}/I_{i+1}$ by
\[
\psi_i(a + I_i) = \phi_i(a) + I_{i+1}.
\]
One can then check that $A/I = \lim_{i \to \infty}(A_i/I_i, \psi_i)$.  It is well 
known that 
\[
A_i/I_i =  \bigoplus_{l=1}^{m_i} p_{i,l}(\mathrm{C}(Y_{i,l}) \otimes \mathcal{K})p_{i,l}
\]
for closed subspaces $Y_{i,l} \subseteq X_{i,l}$, $i \in \mathbb{N}$, $1 \leq l \leq m_i$.
Since $\mathrm{dim}(Y_{i,l}) \leq \mathrm{dim}(X_{i,l})$, we have
\[
\lim_{i \to \infty} \mathrm{max}_{1 \leq l \leq m_i}\left\{\frac{\mathrm{dim}(Y_{i,l})}
{\mathrm{rank}(p_{i,l})}\right\} \leq \mathrm{drr}(A) +\epsilon.
\]
Since $\epsilon$ was arbitrary, we conclude that $\mathrm{drr}(A/I) \leq \mathrm{drr}(A)$.

For (ii), we clearly have $\mathrm{drr}(A \oplus B) \leq \mathrm{max} \{\mathrm{drr}(A),\mathrm{drr}(B)\}$.
$A$ and $B$ are ideals of $A \oplus B$, so we may use (i) to obtain the reverse inequality.  

(iii) is straightforward.

For (iv),
fix decompositions $A = \lim_{i \to \infty}(A_i,\phi_i)$ and $B = \lim_{j \to \infty}(B_j,\psi_j)$,
where
\[
A_i = \bigoplus_{l=1}^{m_i} p_{i,l} (\mathrm{C}(X_{i,l}) \otimes \mathcal{K}) p_{i,l}, \ \ \mathrm{and} \ \
B_j = \bigoplus_{s=1}^{n_i} q_{j,s} (\mathrm{C}(Y_{j,s}) \otimes \mathcal{K}) q_{j,s}.
\]
Assume, as we may, that
\[
\lim_{i \to \infty} \ \mathrm{max}_{1 \leq l \leq m_i}\left\{\frac{\mathrm{dim}(X_{i,l})}{\mathrm{rank(p_{i,l})}}\right\} = \mathrm{drr}(A) +\epsilon_1
\]
and
\[
\lim_{j \to \infty} \ \mathrm{max}_{1 \leq s \leq n_i}\left\{\frac{\mathrm{dim}(Y_{j,s})}{\mathrm{rank(q_{j,s})}}\right\} = \mathrm{drr}(B) + \epsilon_2,
\] 
for some $\epsilon_1,\epsilon_2>0$.
$A \otimes B$ is the limit of the inductive system $(A_i \otimes B_i, \phi_i \otimes \psi_i)$, and 
\[
A_i \otimes B_i = \bigoplus_{l,s} (p_{i,l} \otimes q_{j,s})(\mathrm{C}(X_{i,l} \times
Y_{j,s}) \otimes \mathcal{K})(p_{i,l} \otimes q_{j,s}).
\]
We have the inequalities
\begin{eqnarray}
\mathrm{drr}(A \otimes B) & \leq & \lim_{i \to \infty} \mathrm{max}_{l,s} \left\{ \frac{\mathrm{dim}(X_{i,l}) +
\mathrm{dim}(Y_{i,s})}{\mathrm{rank}(p_{i,l})\mathrm{rank}(q_{i,s})} \right\} \\
& \leq & \lim_{i \to \infty} \mathrm{max}_{l,s} \left\{ \frac{\mathrm{dim}(X_{i,l})}{\mathrm{rank}(p_{i,l})\mathrm{rank}(q_{i,s})}
+  \frac{\mathrm{dim}(Y_{i,s})}{\mathrm{rank}(p_{i,l})\mathrm{rank}(q_{j,s})}\right\} \\
& \leq & \lim_{i \to \infty} \mathrm{max}_{l,s} \left\{ \frac{\mathrm{drr}(A) + \epsilon_1}{\mathrm{rank}(q_{i,s})}
+  \frac{\mathrm{drr}(B) + \epsilon_2}{\mathrm{rank}(p_{i,l})}\right\}.
\end{eqnarray}
Since $A$ and $B$ are simple, we have 
\[
\mathrm{rank}(q_{i,s}), \mathrm{rank}(p_{i,l}) \stackrel{i \to \infty}{\longrightarrow} \infty.  
\]
It follows that the right hand side of (6) is equal to zero, whence $A \otimes B$ has slow dimension growth
and $\mathrm{drr}(A \otimes B)=0$.
\end{proof}

We suspect that equality holds in (iii) above, but it is unclear how a proof might proceed;
a decomposition for $A \otimes \mathrm{M}_k$ need not respect the tensor product structure,
and so does not give rise to an obvious decomposition of $A$.  An inductive limit of 
AH algebras is only approximated locally by semi-homogeneous algebras, and the latter
condition is strictly weaker than approximate homogeneity (\cite{DE}).  Thus, it does
not make sense to investigate the behaviour of the dimension-rank ratio for inductive
limits of AH algebras.

Our first theorem shows that the dimension-rank ratio behaves as one would like for
semi-homogeneous $C^*$-algebras.  Note that the spectrum of $B$ in 
the proposition below need not be a CW-complex, and need not be of finite covering 
dimension.

\begin{thms}\label{drrsemihom}
Let $B = \oplus_{j=1}^n B_j$ be a direct sum of homogeneous $C^*$-algebras
\[
B_j = p_j(\mathrm{C}(X_j) \otimes \mathcal{K})p_j,
\]
Then, 
\[
\mathrm{drr}(B) = \mathrm{max}_{1 \leq j \leq n} \ \left\{ \frac{\mathrm{dim}(X_j)}{\mathrm{rank}(p_j)} \right\}.
\]
\end{thms}

\begin{proof}
By part (i) of Proposition \ref{drrprop}, it will be enough to establish the theorem for $n=1$
and $X_1$ connected.  

Suppose first that $B = p(\mathrm{C}(X) \otimes \mathcal{K})p$ for some connected
compact Hausdorff space $X$ of finite covering dimension.
Clearly,
\[
\mathrm{drr}(B) \leq \frac{\mathrm{dim}(X)}{\mathrm{rank}(p)},
\]
since we may write $B = \lim_{i \to \infty}(B,\mathbf{id}_B)$.
Let $B = \lim_{i \to \infty}(A_i,\phi_i)$ be an injective decomposition for $B$, 
where the $A_i$ and $\phi_i$ are as in
(\ref{stddecomp}) and (\ref{stdAi}).  Let $\mathrm{dr}(\bullet)$ denote the decomposition
rank of a nuclear $C^*$-algebra.  In \cite{Wi3} it is proved that
\[
\mathrm{dr}\left( p(\mathrm{C}(X) \otimes \mathcal{K})p \right) = \mathrm{dim}(X)
\]
whenever $X$ is a compact.  Section 3 of  \cite{KW} shows that 
\[
\mathrm{dr}(C \oplus D) = \mathrm{max} \{\mathrm{dr}(C),\mathrm{dr}(D)\}
\]
for any nuclear $C$ and $D$.  It follows that 
$\mathrm{dr}(B) = \mathrm{dim}(X)$, and that
\[
\mathrm{dr}(A_i) = \mathrm{max}_{1 \leq l \leq m_i} \ \left\{
\mathrm{dim}(X_{i,l}) \right\}.
\]
If $\mathrm{dr}(A_i) \leq n$ for every $i \in \mathbb{N}$, then $\mathrm{dr}(B) \leq n$,
again by section 3 of \cite{KW}.  By dropping terms from the inductive sequence for $B$, we
may assume that $\mathrm{dr}(A_i) = \mathrm{dr}(B)$ for every $i \in \mathbb{N}$.
In other words there exists, for each $i \in \mathbb{N}$, an $1 \leq l_i \leq m_i$
such that 
\[
\mathrm{dim}(X_{i,l_i}) = \mathrm{dim}(X).
\]

If 
\[
\mathrm{rank}(p_{i,l}) > \mathrm{max}_{1 \leq j \leq n} \{\mathrm{rank}(p_j)\},
\]
then the canonical map from $p_{i,l}(\mathrm{C}(X_{i,l}) \otimes \mathcal{K})p_{i,l}$ to
$B$ must be zero, contradicting the injectivity of the decomposition.  Thus, 
\[
\mathrm{rank}(p_{i,l}) \leq  \mathrm{max}_{1 \leq j \leq n} \{\mathrm{rank}(p_j)\}, \ 
\forall i \in \mathbb{N}, \ \forall 1 \leq l \leq m_i.
\]
Suppose that 
\[
\mathrm{max}_{1 \leq l \leq m_i} \ \left\{ \frac{\mathrm{dim}(X_{i,l})}{\mathrm{rank}(p_{i,l})}
\right\} < \frac{\mathrm{dim}(X)}{\mathrm{rank}(p)}.
\]
Then,
\[
\frac{\mathrm{dim}(X_{i,l_i})}{\mathrm{rank}(p_{i,l_i})} <
 \frac{\mathrm{dim}(X)}{\mathrm{rank}(p)},
\]
which, since $\mathrm{dim}(X_{i,l_i}) = \mathrm{dim}(X)$, implies that
\[
\frac{1}{\mathrm{rank}(p_{i,l_i})} < \frac{1}{\mathrm{rank}(p)}.
\]
But this implies that $\mathrm{rank}(p_{i,l_i}) > \mathrm{rank}(p)$, a contradiction.
It follows that
\[
\mathrm{drr}(p(\mathrm{C}(X) \otimes \mathcal{K})p) = \frac{\mathrm{dim}(X)}{\mathrm{rank}(p)}.
\]

If $X$ is infinite-dimensional, then the decomposition rank argument from the second
paragraph of the proof allows us to assume that for each $i \in \mathbb{N}$, there is $1 \leq l_i \leq m_i$
such that 
\[
\mathrm{dim}(X_{i,l_i}) \geq i,
\]
and that the partial map from $p_{i,l_i} (\mathrm{C}(X_{i,l_i}) \otimes \mathcal{K})p_{i,l_i}$ is 
not zero.
On the other hand, rank considerations show that there is $M > 0$ such that
$\mathrm{rank}(p_{i,j}) < M$ whenever the partial map from $p_{i,j} (\mathrm{C}(X_{i,j}) \otimes \mathcal{K})p_{i,j}$
is not zero.  Thus, 
\[
\lim_{i \to \infty} \ \mathrm{max}_{1 \leq l \leq m_i} \ \left\{ 
\frac{\mathrm{dim}(X_{i,l})}{\mathrm{rank}(p_{i,l})}\right\} = \infty
\]
for every decomposition, and $\mathrm{drr}(B)=\infty$, as desired.
\end{proof}

\begin{cors}
Let $A = \lim_{i \to \infty}(A_i,\phi_i)$ be an unital AH algebra, where each
$A_i$ is semi-homogeneous.  Then, $\mathrm{drr}(A) \leq \liminf_{i \to \infty} 
\mathrm{drr}(A_i)$.
\end{cors}

\begin{proof}

There is a sequence  $(n_k)_{k=1}^{\infty}$ of natural numbers such that 
\[
\lim_{k \to \infty} \mathrm{drr}(A_{n_k}) = \liminf_{i \to \infty} \mathrm{drr}(A_i)
\]
in the extended reals, and $A = \lim_{k \to \infty}(A_{n_k},\phi_{n_k})$.  Assuming 
the notation from (\ref{stdAi}) for the $A_{n_k}$s, we have
\[
\mathrm{drr}(A_{n_k}) =  \mathrm{max}_{1 \leq l \leq m_{n_k}} \left\{ 
\frac{\mathrm{dim}(X_{n_k,l})}{\mathrm{rank}(p_{n_k,l})} \right\} 
\stackrel{k \to \infty}{\longrightarrow} \liminf_{i \to \infty} \mathrm{drr}(A_i).
\]
This gives $\mathrm{drr}(A) \leq \liminf_{i \to \infty} \mathrm{drr}(A_i)$ by
definition.
\end{proof}

We conclude this section by noting a connection between the dimension-rank ratio
and Rieffel's stable rank for $C^*$-algebras (\cite{Ri1}).  Let $\mathrm{sr}(A)$ denote
the stable rank of a $C^*$-algebra $A$, and let $\lceil x \rceil$ (resp. $\lfloor x
\rfloor$) denote the least (resp. greatest) integer greater (resp. less) than $x \in \mathbb{R}$.
Consider the following formula, established by Nistor in \cite{N}:
\begin{equation}\label{homsr}
\mathrm{sr}( p(\mathrm{C}(X) \otimes \mathcal{K})p) = \left\lceil \frac {\lfloor \mathrm{dim}(X)/2 \rfloor}
{\mathrm{rank}(p)} \right\rceil + 1
\end{equation}
whenever $X$ is a compact Hausdorff space and $p \in \mathrm{C}(X) \otimes \mathcal{K}$ is
a projection of constant rank.  Clearly,
the right hand side is all but equal to $2 \mathrm{drr}( p(\mathrm{C}(X) \otimes \mathcal{K})p)$,
with any difference owing to the fact that the dimension-rank ratio need not be an integer.
This observation leads to:

\begin{props}
Let $A$ be an unital AH algebra.  
Then, 
\[
\mathrm{drr}(A) \geq \frac{\mathrm{sr}(A)}{2} -1.
\] 
\end{props}

\begin{proof}
The proposition is trivial if $\mathrm{sr}(A)=1,2$.

Suppose that $\mathrm{sr}(A) < \infty$.
Theorem 5.1 of \cite{N} states that if $A = \lim_i(A_i,\phi_i)$ is an inductive limit algebra
where $\mathrm{sr}(A_i) \leq n, \ \forall i \in \mathbb{N}$, then $\mathrm{sr}(A) \leq n$.
Thus, we may assume that regardless of the decomposition $A = \lim_{i \to \infty}(A_i,\phi_i)$, 
one has $\mathrm{sr}(A_i) \geq \mathrm{sr}(A)$.  If the $A_i$ are direct sums of homogeneous
building blocks as in equation (\ref{stdAi}), then by (\ref{homsr}) above we have 
\[
\left\lceil \frac {\lfloor \mathrm{dim}(X_{i,l})/2 \rfloor}{\mathrm{rank}(p_{i,l})} \right\rceil + 1 \geq \mathrm{sr}(A)
\]
for some $1 \leq l \leq m_i$.  Straightforward calculation yields $\mathrm{drr}(A_i) \geq (\mathrm{sr}(A)-2)/2$,
so that $\limsup_{i \to \infty} \mathrm{drr}(A_i) \geq (\mathrm{sr}(A)-2)/2$.  Since the decompostion
of $A$ was arbitrary, we conclude that $\mathrm{drr}(A) \geq (\mathrm{sr}(A)-2)/2$.  

The case of $\mathrm{sr}(A) = \infty$ is similar. 
\end{proof}

\section{Ordered $\mathrm{K}$-theory}

In this section we establish connections between the dimension-rank ratio and
the ordered $\mathrm{K}$-theory of AH algebras.  We
examine first the case of a homogeneous $C^*$-algebra with spectrum a
CW-complex of finite dimension. 
\begin{thms}[Husemoller, Theorems 1.2 and 1.5, Chapter 8, \cite{H}]\label{kstab}
Let $X$ be an $n$-dimensional CW-complex, and let $\gamma, \omega$ be complex
vector bundles over $X$. 
\begin{enumerate}
\item[(i)] If $\gamma$ and $\omega$ are stably isomorphic and the fibre dimension
of $\gamma$ is greater than or equal to $ \lceil n/2 \rceil$, then $\gamma$ and 
$\omega$ are isomorphic.
\item[(ii)] If the fibre dimension of $\gamma$ exceeds that of $\omega$ by an 
amount greater than or equal to $\lceil n/2 \rceil$, then $\omega$ is
isomorphic to a sub-bundle of $\gamma$.
\end{enumerate}
\end{thms}

\noindent
Making the identifications
\[
\mathrm{K}_0(p(\mathrm{C}(X) \otimes \mathcal{K})p) \equiv \mathrm{K}_0(\mathrm{C}(X))
\equiv \mathrm{K}^0(X),
\]
we recast Theorem \ref{kstab} in terms of $\mathrm{K}$-theory (this is standard fare).  Let $p,r \in 
\mathrm{M}_{\infty}(\mathrm{C}(X))$ be projections, and let $[p],[r]$ denote 
their $\mathrm{K}_0$-classes.  Then parts (i) and (ii) of Theorem \ref{kstab} are
equivalent to the following two statements, respectively:
\begin{enumerate}
\item[(i)] if $[p]=[r]$ and $\mathrm{rank}(p) \geq \lceil \mathrm{dim}(X)/2 \rceil$, then
$p$ and $r$ are Murray-von Neumann equivalent;
\item[(ii)] if $\mathrm{rank}(p) - \mathrm{rank}(r) \geq \lceil \mathrm{dim}(X)/2 \rceil$,
then $r$ is Murray-von Neumann equivalent to a subprojection of $p$ ($r \prec p$).  In particular
$[p]-[r] \in \mathrm{K}_0(\mathrm{C}(X))^+$.
\end{enumerate}

Let $A$ be an unital stably finite
$C^*$-algebra, and let $\mathrm{QT}(A)$ denote the compact convex set 
of normalised quasi-traces on $A$. 
(A deep theorem of Haagerup (\cite{Ha}) asserts that every quasi-trace on an
unital and exact $C^*$-algebra $A$ is a trace.  Thus, when
$A$ is exact, unital, and stably finite, we identify $\mathrm{QT}(A)$ with the space $\mathrm{T}(A)$
of normalised traces on $A$.)  We recall three familiar concepts in the $\mathrm{K}$-theory of 
$C^*$-algebras:
\begin{enumerate}
\item[(i)]
If projections $p,q \in \mathrm{M}_{\infty}(A)$ are Murray-von Neumann equivalent whenever
$[p] = [q] \in \mathrm{K}_0A$, then $A$ is said to have {\it cancellation of projections} 
(or simply {\it cancellation}).  
\item[(ii)] If the condition that $\tau(p) < \tau(q)$ for
every $\tau \in \mathrm{QT}(A)$ implies that $p$ is Murray-von Neumann equivalent
to a subprojection of $q$, then we say that $A$ has (FCQ) --- $A$ satisfies 
{\it Blackadar's Second Fundamental Comparability Question}.  
\item[(iii)] If, given elements 
$x_1,x_2,y_1,y_2$ in a partially ordered Abelian group $(G,G^+)$ such that $x_i 
\leq y_j$, $i,j \in \{1,2\}$, there exists $z \in G$ such that $x_i \leq z \leq y_j$, 
$i,j \in \{1,2\}$, then we say that $G$ has the {\it Riesz interpolation property}
(or simply {\it interpolation}).  
\end{enumerate}
Our next definition generalises these notions and another besides.

\begin{dfs}\label{rcancfcq}
Let $A$ be an unital and stably finite $C^*$-algebra, $p,q \in \mathrm{M}_{\infty}(A)$
projections, and $r \geq 0$.  
\begin{enumerate}
\item[(i)] Say that $A$ has $r$-cancellation if $p$ and $q$ are
Murray-von Neumann equivalent whenever $[p]=[q]$ and
\[
\tau(p) = \tau(q) > r, \ \forall \tau \in \mathrm{QT}(A).
\]
\item[(ii)] Say that $A$ has $r$-(FCQ) if $p$ is Murray-von Neumann equivalent
to a subprojection of $q$ whenever
\[
\tau(p) + r < \tau(q), \ \forall \tau \in \mathrm{QT}(A).
\]
\item[(iii)]
Let $(G,G^+,u)$ be a partially ordered Abelian group with distinguished order unit $u$
and state space $S(G)$.  Let $r > 0$.  
Say that $G$ has $r$-interpolation if whenever one has elements
$x_1,x_2,y_1,y_2 \in G$ such that $x_i \leq y_j$, $i,j \in \{1,2\}$, and
\[
s(x_i) + r < s(y_j), \ i,j \in \{1,2\}, \ \forall s \in S(G),
\]
then there exists $z \in G$ such that $x_i \leq z \leq y_j$, $i,j \in \{1,2\}$.
\item[(iv)]
Let $(M,u)$ be a positive ordered semigroup with distinguished strong
order unit $u$ and state space $S(M)$.  Let $r > 0$ and $x,y \in M$.  
Say that $M$ has {\it $r$-strict comparison} if
\[
s(x) + r < s(y), \ \forall s \in S(M), 
\]
implies that $x \leq y$ in $M$.
\end{enumerate} 
\end{dfs}

\noindent
We will prove that the elements of Definition \ref{rcancfcq} are connected naturally
to the dimension-rank ratio.

To prepare the next proposition, recall that a positive ordered semigroup 
$(M,\leq)$ is said to have an {\it algebraic order} if whenever one has
$x,y \in M$ such that $x \leq y$, then there is $z \in M$ such that 
$x+z = y$.  $M$ is said to be {\it cancellative} if whenever one has
elements $x,y,z \in M$ such that $x+z=y+z$, then $x = y$ (cf. \cite{Go2}).

\begin{props}\label{rcomp}
Let $(M,v)$ be a positive ordered semigroup with distinguished strong order unit
$v$.  Suppose that the order on $M$ is algebraic, and that $M$ is cancellative.
Let $G$ be the Grothendieck enveloping group of $M$.  Let $\iota:M \to G$ denote
the Grothendieck map, and put $G^+ = \iota(M)$, $u = \iota(v)$.  Let $S(G)$ denote
the state space of $G$.

Let $r >0$ and $x,y \in G$.  Then,
\[
s(x) + r < s(y),  \ \forall s \in S(G),
\]
implies that $x \leq y$ in $G$ if and only if $(M,v)$ has $r$-strict comparison.
\end{props}

\begin{proof}
Our hypotheses on $M$ imply that $(M,v) \cong (\iota(M),\iota(v)) = (G^+,u)$,
whence $(G^+,u)$ has $r$-strict comparison if and only if $(M,v)$ does.

We may identify $S(G)$ and $S(G^+)$, whence the forward implication follows from
restricting to $G^+$.  (There is a subtle point here:  states on partially
ordered Abelian groups are merely positive homomorphisms into the reals which
take the order unit to $1 \in \mathbb{R}$, whereas states on ordered Abelian 
semigroups are, in addition, order preserving.  We are using the fact that 
$\iota(M) \cong M$ whenever $M$ is algebraically ordered and cancellative to
make our identification of state spaces (\cite{Go2}).  We are grateful to Francesc
Perera for pointing this out to us.)

Now suppose that $(G^+,u)$ has $r$-strict comparison.  Let $x,y \in G$ and
write 
\[
x = x_+ - x_-, \ y = y_+ - y_-, 
\]
where $x_+, x_-,y_+,y_- \in G^+$.  If
\[
s(x) + r < s(y),  \ \forall s \in S(G),
\]
then
\[
s(x_+ + y_-) + r < s(y_+ + x_-),  \ \forall s \in S(G) \equiv S(G^+),
\]
whence $x_+ + y_- \leq y_+ + x_-$ in $G^+$.  It follows that $x \leq y$, as desired.
\end{proof}

In light of the proposition above, we will say that a partially ordered Abelian
group $(G,G^+,u)$ such that $G^+ \cong \iota(G^+)$ has $r$-strict comparison 
whenever $(G^+,u)$ does;  this definition makes sense for the ordered
$\mathrm{K}_0$-group of an unital and stably finite $C^*$-algebra.  

\begin{props}\label{closed} Let $A$ be an unital and stably finite $C^*$-algebra, and $(G,G^+,u)$
a partially ordered Abelian group with distinguished order unit $u$.  Then, the
following sets are closed:
\begin{enumerate}
\item[(i)] $A_1 := \{ r \in \mathbb{R}| \ A \ \mathrm{has} \ r-{\mathrm{cancellation}} \ \}$;
\item[(ii)] $A_2 := \{ r \in \mathbb{R}| \ A \ \mathrm{has} \ r-{\mathrm{(FCQ)}} \ \}$;
\item[(iii)] $A_3 := \{ r \in \mathbb{R}| \ G \ \mathrm{has} \ r-{\mathrm{strict}} \  \mathrm{comparison} \ \}$;
\item[(iv)] $A_4 := \{ r \in \mathbb{R}| \ G \ \mathrm{has} \ r-{\mathrm{interpolation}} \ \}$.
\end{enumerate}
\end{props}

\begin{proof}  For each $i \in \{1,2,3,4\}$ one
has that $s \in A_i$ whenever $s > r$ and $r \in A_i$, so it will suffice to prove that
$\alpha_i := \mathrm{inf}(A_i) \in A_i$.  The proof of each case follows a common thread.

For (i), let there be given projections $p,q \in \mathrm{M}_{\infty}(A)$ such that
\[
[p]=[q], \ \tau(p) = \tau(q) > \alpha_1, \ \forall \tau \in \mathrm{QT}(A).
\] 
The map $\tau \mapsto \tau(p)$ on $\mathrm{QT}(A)$ is continuous and $\mathrm{QT}(A)$
is compact, so this map achieves a minimum value $\delta > \alpha_1$.  Since $\delta \in 
A_1$, we conclude that $p$ and $q$ are Murray-von Neumann equivalent, as desired.

For (ii), let there be given projections $p,q \in \mathrm{M}_{\infty}(A)$ such that
\[
\tau(p) + \alpha_2 < \tau(q), \ \forall \tau \in \mathrm{QT}(A).
\] 
The map $\tau \mapsto \tau(q)-\tau(p)$ is continuous on the compact space $\mathrm{QT}(A)$,
and so achieves a minimum value $\delta > \alpha_2$.  Thus,
\[
\tau(p) + \delta < \tau(q), \ \forall \tau \in \mathrm{QT}(A).
\] 
Since $\delta \in A_2$, the desired conclusion follows.

For (iii), let $x,y \in G^+$ be such that 
\[
s(x) + \alpha_3 < s(y), \ \forall s \in S(G).
\]
The map
\[
s \mapsto s(y) - s(x)
\]
is strictly positive and continuous, and the space $S(G)$ is compact (cf. Proposition 6.2, \cite{Go2}).  
Thus, this map achieves a minimum value $\delta > \alpha_3$.  We now have
\[
s(x) + \delta < s(y), \ \forall s \in S(G).
\]
Since $\delta \in A_3$, the desired conclusion follows.

For (iv), let there be given elements $x_1,x_2,y_1,y_2 \in G$ satisfying
\[
s(x_i) + \alpha_4 < s(y_j), \ i,j \in \{1,2\}, \ \forall s \in S(G).
\]
For each pair $(i,j)$, $i,j \in \{1,2\}$, there exists $r_{i,j} > 0$ such that
\[
s(x_i) + \alpha_4 + r_{i,j} < s(y_j), \ \forall s \in S(G).
\]
Put $\delta = \alpha_4 + \mathrm{min} \{r_{1,1},r_{1,2},r_{2,1},r_{2,2} \}$. 
Now
\[
s(x_i) + \delta < s(y_j), \ i,j \in \{1,2\}, \ \forall s \in S(G),
\]
and $\delta \in A_4$.  We conclude that there is an interpolating element
$z \in G$ such that $x_i \leq z \leq y_j$, $\forall i,j \in \{1,2\}$.
\end{proof} 

Definition \ref{rcancfcq} can be used to summarise
the natural connections between the $\mathrm{K}$-theory of homogeneous $C^*$-algebras
and their dimension-rank ratios.  

\begin{props}\label{rhomprop}
Let $A = p(\mathrm{C}(X) \otimes \mathcal{K})p$, where $X$ is a connected CW-complex of finite dimension.  
Then: 
\begin{enumerate}
\item[(i)] $A$ has $(\mathrm{drr}(A)/2)$-cancellation;
\item[(ii)] $A$ has $(\mathrm{drr}(A)/2)$-(FCQ);
\item[(iii)] $\mathrm{K}_0A$ has $(\mathrm{drr}(A)+1/\mathrm{rank}(p))$-interpolation;
\item[(iv)] $\mathrm{K}_0A^+$ has $(\mathrm{drr}(A)/2)$-strict comparison.
\end{enumerate}
\end{props}

\begin{proof}
(i), (ii), and (iv) are straightforward: combine Definition \ref{rcancfcq} with Theorem \ref{kstab}, (ii).  
We prove (iii), which is slightly more involved.

Let $s$ denote the unique (geometric) state on $\mathrm{K}_0A$, and recall
that for a projection $r \in \mathrm{M}_{\infty}(A)$ we have
\[
s([r]) = \frac{\mathrm{rank}(r)}{\mathrm{rank}(p)}.
\]
For the remainder of the proof, let $r,q \in \mathrm{M}_{\infty}(A)$ be projections.

Assume that we are
given four elements $x_{1},x_{2},y_{1},y_{2} \in \mathrm{K}_{0}A$ such that
$x_{i} \leq y_{j}$, $i,j \in \{1,2\}$, and
\begin{equation}\label{blah}
s(x_{i}) + \mathrm{drr}(A) + 1/\mathrm{rank}(p) < s(y_{j}), \ i,j \in \{1,2\}.
\end{equation}
Every element $x \in \mathrm{K}_{0}A$ can be written as a difference of
$\mathrm{K}_{0}$-classes of projections, say $x = [q]-[r]$.  The
difference $\mathrm{rank}(q) - \mathrm{rank}(r)$ is commonly referred 
to as the {\it virtual dimension} of $x$.  We will let
$\mathrm{rank}(x)$ denote this virtual dimension, thus extending the notion
of rank to all of $\mathrm{K}_{0}A$.  With this notation we have
\[
s(x) = \frac{\mathrm{rank}(x)}{\mathrm{rank}(p)}, \ \forall x \in
\mathrm{K}_{0}A.
\]
We may now rewrite (\ref{blah}) above as
\[
\frac{\mathrm{rank}(x_i)}{\mathrm{rank}(p)} + \frac{\mathrm{dim}(X)}{\mathrm{rank}(p)}
+ \frac{1}{\mathrm{rank}(p)} < \frac{\mathrm{rank}(y_j)}{\mathrm{rank}(p)},
\]
which yields
\[
\mathrm{rank}(y_j) - \mathrm{rank}(x_i) > \mathrm{dim}(X) +1, \ i,j \in \{1,2\}.
\]
Let $z$ be any element of $\mathrm{K}_0A$ such that 
\[
\mathrm{rank}(z) = \mathrm{max} \{\mathrm{rank}(x_1),\mathrm{rank}(x_2)\} + \lceil \mathrm{dim}(X)/2 \rceil.
\]
Then 
\[
\mathrm{rank}(z-x_i), \mathrm{rank}(y_j-z) \geq \lceil \mathrm{dim}(X)/2 \rceil, \ i,j \in \{1,2\},
\]
and $z$ is the desired interpolating element by Theorem \ref{kstab}, (ii).

\end{proof}

\noindent
We shall see below that Proposition \ref{rhomprop} can be generalised to the 
setting of general unital AH algebras, provided that the algebras have 
ordered $\mathrm{K}_0$-groups which admit a unique state.

\begin{exs} {\rm
While part (iii) of Proposition \ref{rhomprop} gives a positive real $r$ 
such that the algebra $A$ as in the hypotheses has $r$-interpolation, it is not immediately
clear that there may be a nonzero lower bound on the set of all such reals.  But be
one there may.  Consider, for any natural number $n > 1$, the $C^*$-algebra 
$A = \mathrm{M}_n(\mathrm{C}(S^{2n}))$.
Clearly, $\mathrm{drr}(A) = 2$.  The ordered $\mathrm{K}_0$-group of $A$ is 
well known:  it is isomorphic as a group to $\mathbb{Z} \oplus \mathbb{Z}$;
the first co-ordinate is generated by the $\mathrm{K}^0$-class $[\theta_1]$ of the 
trivial line bundle $\theta_1$;  the second co-ordinate is generated by the 
difference $[\xi]-[\theta_n]$, where $\xi$ is the bundle corresponding to the 
$n$-dimensional Bott projection and $\theta_n$ is the trivial bundle of fibre
dimension $n$;  the positive cone $\mathrm{K}_0A^+$ is
\[
\{(x,y)| y=0 \ \mathrm{and} \ x \geq 0 \} \cup \{ (x,y)| x \geq n \}.
\]
Put
\[
x_1 = 0 \oplus 0, \ x_2 = 0 \oplus 1, \ y_1 = n \oplus 0, \ y_2 = n \oplus 1.
\]
With the description of $\mathrm{K}_0A^+$ in hand, one checks easily that
\[
x_i \leq y_j, \ i,j \in \{1,2\}.
\]
in $\mathrm{K}_0A$, yet there is no $z \in \mathrm{K}_0A$ which interpolates
these four elements.  The unique geometric state $s$ on $\mathrm{K}_0A$ returns
the rank of a $\mathrm{K}_0$ element divided by $n$, whence 
\[
s(x_i) + r < s(y_j), \ \forall r < 1, \ i,j \in \{1,2\}.
\]
Thus, $\mathrm{K}_0A$ does not have $r$-interpolation
for any $r < \mathrm{drr}(A)/2$. }
\end{exs}

\begin{lms}\label{onestate}
Let $A \cong \lim_{i \to \infty}(A_i,\phi_i)$ be an unital AH algebra, where each $A_i$
is homogeneous with connected spectrum.  Then, $(\mathrm{K}_0(A),\mathrm{K}_0(A)^+,[1_A])$ is a simple partially
ordered Abelian group admitting a unique state.
\end{lms}

\begin{proof}
For each $i \in \mathbb{N}$ write
\[
A_i = p_i(\mathrm{C}(X_i) \otimes \mathcal{K})p_i,
\]
where $X_i$ is a compact connected Hausdorff space,
and $p_i \in \mathrm{C}(X_i) \otimes \mathcal{K}$ is a projection.  As noted following
Definition \ref{drr}, the $X_i$ may be assumed to have finite covering dimension (cf. \cite{Bl2},
\cite{Go1}).  $(\mathrm{K}_0(A), \mathrm{K}_0(A)^+)$ is a partially ordered Abelian group
for every stably finite $A$ (cf. \cite[Chapter 6, Section 3]{Bl1}).  

There is a unique (geometric) state on $\mathrm{K}_0 A_i$ which returns the normalised rank of
a projection corresponding to a positive $\mathrm{K}_0$-class, and is extended to
all of $\mathrm{K}_0A_i$ by linearity.  By Proposition 6.14 of \cite{Go2}, 
$S(\mathrm{K}_0 A)$ is the inverse limit of the $S(\mathrm{K}_0 A_i)$s,
whence $\mathrm{K}_0 A$ admits a unique state.

It remains to prove that $(\mathrm{K}_0(A),\mathrm{K}_0(A)^+)$ is a simple ordered group, i.e., that every
non-zero positive element is an order unit.  It will suffice to prove that each $(\mathrm{K}_0(A_i),\mathrm{K}_0(A_i)^+)$
is a simple ordered group.  Each element of 
$\mathrm{K}_0(A_i) = \mathrm{K}^0(X_i)$ corresponds to a difference $x = [q] - [p]$, where
$q,p \in \mathrm{M}_{\infty}(A_i)$ are projections.  Let $y = [e]-[f] \in \mathrm{K}_0A_i$,
where $e,f \in \mathrm{M}_{\infty}(A_i)$ are projections.  If $x$ is positive, then
$\mathrm{rank}(q) > \mathrm{rank}(p)$.  In particular, there exists $n \in \mathbb{N}$ such that
\[
\mathrm{rank}(nx) - (\mathrm{rank}([e]) - \mathrm{rank}([f])) \geq \lceil \mathrm{dim}(X_i)/2 \rceil,
\]
so $ny \geq y$ by Theorem \ref{kstab}, (ii), and $x$ is an order unit.
\end{proof}

We will need the following result to prove our next lemma:

\begin{thms}[Goodearl, Proposition 4.16, \cite{Go2}]\label{infsup}
Let $(G,G^+,u)$ be a non-zero partially ordered Abelian group with distinguished order unit.  If $G$
admits a unique state $s$, then for any $x \in G^+$ one has
\begin{eqnarray*}
s(x) & = & \mathrm{inf} \{l/n|l,n \in \mathbb{N} \ \mathrm{and} \ nx \leq lu \} \\
     & = & \mathrm{sup} \{k/m|k \in \mathbb{Z}^+, m \in \mathbb{N}, \ \mathrm{and} \ ku \leq mx \}
\end{eqnarray*}
\end{thms}

\begin{lms}\label{evenstate}
Let $A$ be an unital AH algebra, and suppose that $\mathrm{K}_0(A)$ admits
a unique state.  Let there be given a decomposition of $A$ as in equations
(\ref{stddecomp}) and (\ref{stdAi}) and a tolerance $\epsilon > 0$.  Then,
for any $x \in \mathrm{K}_0(A)^+$ there exists $j \in \mathbb{N}$ such that
\begin{enumerate}
\item[(i)] $x$ has a pre-image $x_j \in \mathrm{K}_0(A_j)$;
\item[(ii)] if
\[
A_j = \bigoplus_{l=1}^{m_j} p_{j,l} (\mathrm{C}(X_{j,l}) \otimes \mathcal{K}) p_{j,l},
\]
and $s_l$ denotes the state on $\mathrm{K}_0(A_j)$ which is equal to the (unique) geometric state
$g_l$ on $\mathrm{K}_0(p_{j,l} (\mathrm{C}(X_{j,l}) \otimes \mathcal{K}) p_{j,l})$ and zero on the other
direct summands of $A_j$, then
\[
|s_l(x) -s(x)| < \epsilon, \ 1 \leq l \leq m_j.
\]
\end{enumerate}
\end{lms}

\begin{proof}
By truncating the given inductive sequence for $A$, we may assume that $x$ has
a pre-image in every $A_i$, $i \in \mathbb{N}$.

Using Theorem \ref{infsup}, find non-negative integers $r,n,k,m$ such that
\[
r/n - s(x) < \epsilon/2, \ \ s(x) - k/m < \epsilon/2,
\]
$nx \leq r[1_{A}]$, and $k[1_{A}] \leq mx$ inside $\mathrm{K}_0(A)$.  The last
two inequalities must hold already in some $A_j$, and, since $\mathrm{K}_0(A_j)$
has the direct sum order coming from the summands $\mathrm{K}_0(p_{j,l} 
(\mathrm{C}(X_{j,l}) \otimes \mathcal{K}) p_{j,l})$, they will still hold upon
restricting to any such summand.
Let $x_l$ denote the restriction of $x$ to $\mathrm{K}_0(p_{j,l} 
(\mathrm{C}(X_{j,l}) \otimes \mathcal{K}) p_{j,l})$.  We have
\[
nx_l \leq r[p_{j,l}], \ \ k[p_{j,l}] \leq mx_l, \ 1 \leq l \leq m_j.
\]
Since the geometric state $g_l$ on $\mathrm{K}_0(p_{j,l} 
\mathrm{M}_{k_{j,l}}(\mathrm{C}(X_{j,l})) p_{j,l})$ preserves order,
we conclude that  
\[
k/m \leq g_l(x_l) \leq r/n.
\]
Since $g_l(x_l) = s_l(x)$, the lemma follows.
\end{proof}

\begin{thms}\label{drrstab}
Let $A$ be an unital AH algebra with $\mathrm{drr}(A) < \infty$, and suppose 
that $\mathrm{K}_0(A)$ admits a unique state $s$.  Then:
\begin{enumerate}
\item[(i)] $A$ has $(\mathrm{drr}(A)/2)$-cancellation;
\item[(ii)] $A$ has $(\mathrm{drr}(A)/2)$-(FCQ);
\item[(iii)] $\mathrm{K}_0A$ has $(\mathrm{drr}(A)/2)$-strict comparison.
\end{enumerate}
If, in addition, $A$ is simple, then
\begin{enumerate}
\item[(iv)] $\mathrm{K}_0(A)$ has $(\mathrm{drr}(A))$-interpolation.
\end{enumerate}
\end{thms}

\begin{proof}
We prove that $A$ has $(\mathrm{drr}(A)/2+\epsilon)$-cancellation,
$(\mathrm{drr}(A)/2+\epsilon)$-(FCQ), $(\mathrm{drr}(A)/2+\epsilon)$-strict 
comparison, and $(\mathrm{drr}(A) + \epsilon)$-interpolation for every $\epsilon>0$;  the theorem then
follows from Proposition \ref{closed}.  Let $\epsilon > 0$ be given.

For (i), let there be given projections $p,q \in \mathrm{M}_{\infty}(A)$ such that
$[p]=[q]$ and 
\[
\tau(p) = \tau(q) = s([p]) = s([q]) > \mathrm{drr}(A)/2 + \epsilon, \ \forall \tau \in \mathrm{T}(A).
\]
Fix a decomposition $A=\lim_{i \to \infty}(A_i,\phi_i)$ where
\[
A_{i} =  \bigoplus_{l=1}^{m_{i}} p_{i,l} (\mathrm{C}(X_{i,l}) \otimes \mathcal{K}) p_{i,l}
\]
and 
\[
\mathrm{max}_{1 \leq l \leq m_i} \left\{ \frac{\mathrm{dim}(X_{j,l})}{\mathrm{rank}(p_{j,l})} 
\right\} \leq \mathrm{drr} + \epsilon/2, \ \forall i \in \mathbb{N}.
\]
Use Lemma \ref{evenstate}
to find $j \in \mathbb{N}$ such that $p$ and $q$ have pre-images at the level of $\mathrm{K}_0$ 
(which are projections) $\tilde{p}$ and
$\tilde{q}$, respectively, in $\mathrm{M}_{\infty}(A_j)$ with the properties that
$[\tilde{p}] = [\tilde{q}]$, $s([\tilde{p}]) = s([p])$, and 
\[
|s_l([\tilde{p}])-s([\tilde{p}])| = |s_l([\tilde{q}])-s([\tilde{q}])| < \frac{s(p) - \mathrm{drr}(A)}{4}, \ 
1 \leq l \leq m_j.
\]
Since $s_l([\tilde{p}])$ represents the normalised rank of $\tilde{p}$ restricted to
the direct summand $p_{j,l} (\mathrm{C}(X_{j,l}) \otimes \mathcal{K}) p_{j,l}$ of $A_j$, we
conclude that this restriction is in
the stable range of $\mathrm{K}_0(p_{j,l} (\mathrm{C}(X_{j,l}) \otimes \mathcal{K}) p_{j,l})$
(and similarly for the restriction of $\tilde{q}$).
Thus, the said restrictions, having the same class in $\mathrm{K}_0$, 
are Murray-von Neumann equivalent by Theorem \ref{kstab}.  
It follows that $\tilde{p}$ and $\tilde{q}$ are Murray-von Neumann
equivalent, whence so are $p$ and $q$.  This shows that $A$ has 
$(\mathrm{drr}(A)/2+\epsilon)$-cancellation.  Since $\epsilon$ was arbitrary,
this proves (i).

For (ii), (iii), and (iv) we will retain the decomposition of $A$ from the proof of (i);
for (ii) and (iii) we will retain as the pre-images $\tilde{p}$ and $\tilde{q}$ of $p$ 
and $q$ above, with the property that
\[
|s_l([\tilde{p}])-s([\tilde{p}])|, \ |s_l([\tilde{q}])-s([\tilde{q}])| < \frac{s(p) - \mathrm{drr}(A)}{4}, \ 
1 \leq l \leq m_j.
\]

For (ii), let there be given projections $p,q \in \mathrm{M}_{\infty}(A)$ such that 
\[
\tau(p) + \mathrm{drr}(A)/2 + \epsilon < \tau(q), \ \forall \tau \in \mathrm{T}(A).
\]
Since $\mathrm{K}_0A$ has a unique state $s$, the statement above is equivalent to
\[
s([p]) + \mathrm{drr}(A)/2 + \epsilon < s([q]).
\]
Find pre-images $\tilde{p}$ and $\tilde{q}$ as before.  Then, the virtual dimension
of the restriction of $[\tilde{q}]-[\tilde{p}]$ to a direct summand  $p_{j,l} 
(\mathrm{C}(X_{j,l}) \otimes \mathcal{K}) p_{j,l}$ of $A_j$ is in the stable range
of $\mathrm{K}_0(p_{j,l} (\mathrm{C}(X_{j,l}) \otimes \mathcal{K}) p_{j,l})$, whence
the said restriction is positive.  The direct sum of these restrictions, namely, 
$[\tilde{q}]-[\tilde{p}]$ itself, is then positive.  Write $[\tilde{q}] = [\tilde{p}] 
+ [r]$ for some projection $r \in \mathrm{M}_{\infty}(A_j)$.  Since 
\[
s([\tilde{p}]) + s([r]) = s([\tilde{q}]) = s([q]) > \mathrm{drr}(A)/2 + \epsilon,
\]
we conclude by (i) that $\tilde{p} \oplus r$ and $\tilde{q}$ are Murray-von Neumann
equivalent.  It follows that $\tilde{p}$ is equivalent to a subprojection of $\tilde{q}$,
and similarly for $p$ and $q$.  This proves that $A$ has $(\mathrm{drr}(A)/2+\epsilon)$-(FCQ), 
and so proves (ii).

$\mathrm{K}_0A$ has $(\mathrm{drr}(A)/2 +\epsilon)$-strict comparison if and only if
the same is true of the semigroup $(\mathrm{K}_0A^+,[1_A])$.  The latter condition is equivalent
to the statement that for $[p],[q] \in \mathrm{K}_0A^+$ such that
\[
s([p]) + \mathrm{drr}(A)/2 + \epsilon < s([q]), 
\]
one has $[p] \leq [q]$.  This, in turn, follows from (i), proving (iii).

For (iv), we must prove that for any $x_1,x_2,y_1,y_2 \in \mathrm{K}_0A$ such
that
\[
x_i \leq y_j, \ s(x_i) + \mathrm{drr}(A) +\epsilon < s(y_j), \ i,j \in \{1,2\},
\]
there exists $z \in \mathrm{K}_0A$ such that $x_i \leq z \leq y_j$, $i,j \in \{1,2\}$.
We may assume that $x_1 = 0$ and put $x_2 = x$, for convenience --- (iv) then follows by translating $z$.

Fix projections $p_{y_1},p_{y_2},p_x^+,p_x^- \in \mathrm{M}_{\infty}(A)$ such that
\[
y_1 = [p_{y_1}], \ y_2 = [p_{y_2}]; \ x = [p_x^+]-[p_x^-].
\]
Find, as in the proof of (i), some $j \in \mathbb{N}$ such that $p_{y_1},p_{y_2},p_x^+$, and $p_x^-$
have pre-images (at the level of $\mathrm{K}_0$) $\tilde{p}_{y_1}, \tilde{p}_{y_2}, \tilde{p}_x^+$, 
and $\tilde{p}_x^-$ (all projections), respectively, in $\mathrm{M}_{\infty}(A_j)$, with the property that
\[
|s_l([q])-s([q])| < \frac{\epsilon}{4}, \ \forall q \in \{p_{y_1},p_{y_2},p_x^+,p_x^-\}, \ 1 \leq l \leq m_j.
\]
We may assume, by the simplicity of $A$, that $j$ has also been chosen large enough to ensure that
$1/\mathrm{rank}(p_{j,l}) \ll \epsilon/4$, $1 \leq l \leq m_j$.  

Fix a summand $A_{j,l} = p_{j,l}(\mathrm{C}(X_{j,l}) \otimes \mathcal{K})p_{j,l}$ of $A_j$.  This, by
Proposition \ref{rhomprop}, (iii), has 
\[
\frac{\mathrm{dim}(X_{j,l})}{\mathrm{rank}(p_{j,l})} + \frac{1}{\mathrm{rank}(p_{j,l})}
\leq \mathrm{drr}(A) + \epsilon/2 + \epsilon/4 = \mathrm{drr}(A) + 3\epsilon/4
\]
interpolation.  The restrictions of $\tilde{p}_{y_1}, \tilde{p}_{y_2}, \tilde{p}_x^+$, and $\tilde{p}_x^-$
to $A_{j,l}$ are such that:
\begin{eqnarray*}
\mathrm{drr}(A) + 3\epsilon/4 & < & s_l([\tilde{p}_{y_k}|_{A_{j,l}}]), \ k \in \{1,2\}; \\
s_l([\tilde{p}_x^+|_{A_{j,l}}]-[\tilde{p}_x^-|_{A_{j,l}}]) + \mathrm{drr}(A) + 3\epsilon/4 & < & s_l([\tilde{p}_{y_k}|_{A_{j,l}}]), \ k \in \{1,2\}.
\end{eqnarray*}
It follows that there exists $z_l \in \mathrm{K}_0A_{j,l}$ such that 
\[
0, \ [\tilde{p}_x^+|_{A_{j,l}}]-[\tilde{p}_x^-|_{A_{j,l}}] \leq z_l \leq [\tilde{p}_{y_1}|_{A_{j,l}}], \ [\tilde{p}_{y_2}|_{A_{j,l}}].
\]
Thus,
\[
0, \ [\tilde{p}_x^+]-[\tilde{p}_x^-] \leq \oplus_{l=1}^{m_j} z_l \leq [\tilde{p}_{y_1}], \ [\tilde{p}_{y_2}]
\]
in $\mathrm{K}_0A_j$, and, upon taking images in $\mathrm{K}_0A$ and setting $z = \mathrm{K}_0(\phi_{j\infty})(\oplus_{l=1}^{m_j} z_l)$,
\[
0, \ x \leq z \leq y_1, \ y_2,
\]
as desired.
\end{proof}

\section{A classification result}

Clearly, slow dimension growth implies $\mathrm{drr} = 0$.  This begs the obvious
question:

\begin{qus}\label{drrzerosdg?}
Does $\mathrm{drr}(A)=0$ imply that $A$ has slow dimension growth for every
simple unital AH algebra $A$?
\end{qus}

\noindent
The next theorem and corollary provide a positive answer to Question \ref{drrzerosdg?}
in the case of simple algebras with real rank zero.  It is plausible that this positive
answer will extend to simple algebras of real rank one, too.
Recall that a simple partially ordered Abelian group $(G,G^+)$ is said to be {\it weakly unperforated}
if $mx > 0$ for some $m \in \mathbb{N}$ and $x \in G$ implies that $x > 0$ (Chapter 6, \cite{Bl1}).

\begin{thms}\label{wucanc}
Let $A$ be an unital AH algebra such that $\mathrm{drr}(A) = 0$, and suppose that $\mathrm{K}_0A$ is
a simple ordered group.  Then, $\mathrm{K}_0A$ is weakly unperforated, and $A$ has cancellation.
\end{thms}

\begin{proof}
Suppose that $mx > 0$ for some $m \in \mathbb{N}$ and $x \in G$.  Since $\mathrm{K}_0A$ is 
a simple ordered group, there exists $n \in \mathbb{N}$ such that $nmx > [1_A] \in \mathrm{K}_0A$.
Since $\mathrm{drr}(A) = 0$, we may choose an injective decomposition 
\[
A \cong \lim_{i \to \infty} \left( A_i := \bigoplus_{l=1}^{m_i} p_{i,l}(\mathrm{C}(X_{i,l}) \otimes 
\mathcal{K})p_{i,l}, \phi_i \right)
\]
with the property that for every $i \in \mathbb{N}$,
\[
\mathrm{max}_{1 \leq l \leq m_i} \ \left\{ \frac{\mathrm{dim}(X_{i,l})}{\mathrm{rank}(p_{i,l})} \right\} < \frac{1}{nm}.
\]
Find a pre-image $x_i \in \mathrm{K}_0A_i$ of $x$ such that $nmx_i > [1_{A_i}]$.  Write $x_i = [p]-[q]$ for projections
$p$ and $q$ in $\mathrm{M}_{\infty}(A_i)$.  Let $Sp(\bullet)$ denote spectrum of a $C^*$-algebra.
Upon restricting to any direct summand
$B$ of $A_i$ corresponding to a connected component of $Sp(A_i)$ one has
\[
\mathrm{rank}(p|_B) - \mathrm{rank}(q|_B) > \frac{\mathrm{rank}(1_{B})}{nm} \geq \mathrm{dim}(Sp(B)).
\]
It follows from Theorem \ref{kstab} that
$[p|_B]-[q|_B] \in \mathrm{K}_0B^+$, whence $x_i$ and its image $x \in \mathrm{K}_0A$ are positive.
Thus, $\mathrm{K}_0A$ is weakly unperforated.

Now suppose that we are given projections $p,q \in \mathrm{M}_{\infty}(A)$ such that 
$[p]=[q] \in \mathrm{K}_0A$.  Since $\mathrm{K}_0A$ is a simple ordered group, every positive
element is an order unit.  Hence, there exists some $m \in \mathrm{N}$ such that $m[p] = m[q]
\geq [1_A]$.  Find an injective decomposition for $A$ as above, with the property that
for every $i \in \mathbb{N}$,
\[
\mathrm{max}_{1 \leq l \leq m_i} \ \left\{ \frac{\mathrm{dim}(X_{i,l})}{\mathrm{rank}(p_{i,l})} \right\} < \frac{1}{m}.
\]
Find projections $p_i,q_i \in \mathrm{M}_{\infty}(A_i)$, some $i \in \mathbb{N}$, such
that $p_i$ is a pre-image of $p$, $q_i$ is a pre-image of $q$, $[p_i] = [q_i] \in \mathrm{K}_0A_i$,
and $m[p_i] \geq [1_{A_i}] \in \mathrm{K}_0A_i$.
Now, upon restricting to any direct summand
$B$ of $A_i$ corresponding to a connected component of the spectrum of $A_i$ one has
\[
\mathrm{rank}(p_i|_B), \ \mathrm{rank}(q_i|_B) > \frac{\mathrm{rank}(1_{B})}{m} \geq \mathrm{dim}(Sp(B)).
\]
It follows that $p_i|_B$ and $q_i|_B$ are in the stable range of $\mathrm{K}_0B$, whence they 
are Murray-von Neumann equivalent by Theorem \ref{kstab}.  It follows that $p_i$ and
$q_i$ are Murray-von Neumann equivalent, and so are $p$ and $q$.  Thus, $A$ has cancellation.
\end{proof}

This is the natural point at which to prove the next corollary, but its statement
refers to the almost unperforation of the Cuntz semigroup $W(A)$;  we have yet to 
remind the reader of this notion. As we will have occasion to discuss this notion 
in depth in section 6, we defer our definition until then.

\begin{cors}\label{taf}
Let $A$ be a simple unital AH algebra of real rank zero.  Then, the following are equivalent:
\begin{enumerate}
\item[(i)] $\mathrm{drr}(A) = 0$;
\item[(ii)] $A$ is tracially AF;
\item[(iii)] $A$ has slow dimension growth;
\item[(iv)] $A$ is $\mathcal{Z}$-stable;
\item[(v)] $W(A)$ is almost unperforated and $\mathrm{sr}(A)=1$.
\end{enumerate}
\end{cors}

\begin{proof}
The equivalence of (ii), (iii), (iv), and (v) is Theorem 3.13 of \cite{TW2},
and is the work of many hands, including Marius D\u{a}d\u{a}rlat, George Elliott, Guihua Gong, Huaxin Lin, 
Mikael R{\o}rdam, Wilhelm Winter, and the author.  

If $A$ has slow
dimension growth, then $\mathrm{drr}(A) = 0$ by definition.  Thus, (iii)
implies (i)

We now prove that (i) implies (ii).
It follows from Theorem \ref{wucanc} that $\mathrm{K}_0(A)$ is weakly
unperforated and has cancellation of projections.  Combining this with real rank
zero yields stable rank one for $A$ (Proposition 6.5.2, \cite{Bl1}).  
That $A$ is tracially AF then follows from \cite{Li4}.
\end{proof}

\noindent
Corollary \ref{taf} allows us to view the dimension-rank ratio as a measure of dimension growth which 
extends the existing notion of slow dimension growth.  The condition
$\mathrm{drr}=0$ is a more natural way to view slow dimension growth,
since it has higher analogues in the form of non-zero dimension-rank ratios.
As promised, the dimension-rank ratio is insensitive to differences between
$\mathcal{Z}$-stable algebras, provided that they are simple and of real rank zero.

\section{The range of the dimension-rank ratio}

It is clear from Theorem \ref{drrsemihom} that the dimension-rank ratio may take any finite,
nonnegative, and rational value.  In fact, more is true:

\begin{thms}\label{drrrange}
Let $c \in \mathbb{R}^+ \cup \{\infty\}$.  There exists a simple, unital AH algebra $A_c$
such that $\mathrm{K}_0A_c$ admits a unique state and $\mathrm{drr}(A_c)=c$.  Moreover,
the stable rank of $A_c$ is one. 
\end{thms}

\begin{proof}
We address the extreme cases first.
The case $c=0$ is straightforward: any UHF algebra has $\mathrm{drr} = 0$.
For $c = \infty$, we use an existing example due to Villadsen.
In \cite{V1}, Villadsen constructs several simple unital
AH algebras whose $\mathrm{K}_0$-groups admit a unique state $s$.  One of these,
say $A$, has unbounded perforation in its ordered $\mathrm{K}_0$-group ---
for every $n \in \mathbb{N}$, there is a non-positive element $x_n \in \mathrm{K}_0A$
such that $s(x_n) \geq n$.  No matter 
how one decomposes $A$ as an inductive limit of direct sums of homogeneous $C^*$-algebras ---
as 
\[
A = \lim_{i \to \infty} (A_i := \bigoplus_{l=1}^{m_i} p_{i,l}(\mathrm{C}(X_{i,l}) 
\otimes \mathcal{K}) p_{i,l} ,\phi_i), 
\]
say ---
one will always have $x_n$ arising in the $\mathrm{K}_0$-group of $A_j$ for all $j$ greater
than or equal to some $j_0 \in \mathbb{N}$.  Since $\mathrm{K}_0A$ admits a unique state,
we may apply Lemma \ref{evenstate} to conclude that for any $\epsilon > 0$, there is some 
$j \geq j_0$ with the following property:  the restriction $x_{n,l}$ of $x_n$ to the $\mathrm{K}_0$-group of the direct summand  
$p_{j,l}(\mathrm{C}(X_{j,l}) \otimes \mathcal{K})p_{j,l}$ of $A_j$ satisfies
\[
|s_l(x_{n,l}) - s(x_n)| < \epsilon, \ \forall 1 \leq l \leq m_j.
\]
Since 
\[
s_l(x_{n,l}) = \frac{\mathrm{rank}(x_{n,l})}{\mathrm{rank}(p_{j,l})},
\]
we have 
\[
|\mathrm{rank}(x_{n,l}) - s(x_n) \cdot \mathrm{rank}(p_{j,l}) | < \epsilon \cdot \mathrm{rank}(p_{j,l})
\]
and
\[
\mathrm{rank}(x_{n,l}) \geq (s(x_n) - \epsilon) \cdot \mathrm{rank}(p_{j,l})  \geq (n-\epsilon) \cdot \mathrm{rank}(p_{j,l}).
\]
It follows that from Theorem \ref{kstab}, (ii) (rephrased in $\mathrm{K}$-theoretic terms)
that 
\[
\frac{\mathrm{dim}(X_{j,l})}{\mathrm{rank}(p_{j,l})} > \frac{n-1}{2}, \ \forall 1 \leq l \leq m_j.
\]
Since $n$ was arbitrary, we conclude that no matter the decomposition, 
$\limsup_{i \to \infty} \mathrm{drr}(A_i) = \infty$; $\mathrm{drr}(A) = \infty$
by definition.

Now suppose that $c \in \mathbb{R}^+ \backslash \{0\}$.  We construct $A_c$ by methods similar to those 
of \cite{V1}.  $A_c$ will be the limit of an inductive sequence $(B_i,\phi_i)$,
where 
\[
B_i = \mathrm{M}_{n_i}(\mathrm{C}(X_i)), \ X_i = (\mathrm{S}^2)^{m_1 m_2 \cdots m_i},
\] 
and $n_i,m_i \in \mathbb{N}$ are to be specified.

Choose $m_1$ and $n_1$ so that $m_1/n_1 > c/2$.  We have $X_{i+1} = (X_i)^{m_{i+1}}$ by
construction.  Let
\[
\pi_i^j: (X_i)^{m_{i+1}} \to X_i, \ 1 \leq j \leq m_{i+1},
\]
be the co-ordinate projections.  Define a map $\phi_i:B_i \to B_{i+1}$ by  
\[
\phi_i(f)(x) = \mathrm{diag}\left(f \circ \pi_i^1(x), \ldots,f \circ \pi_i^{m_{i+1}}(x), f(x_i^1), \ldots f(x_i^{s_{i+1}})\right),
\]
where $s_{i+1} \in \mathbb{N}$ and the $x_i^1,\ldots,x_i^{s_{i+1}} \in X_i$ are to be specified.  
Suppose that for $i \leq k$ we have chosen the parameters in our construction inductively so that
\begin{equation}\label{closetodrr}
\frac{c}{2} < \frac{m_1 m_2 \cdots m_k}{n_k} < \frac{c}{2} + \frac{1}{2^k}.
\end{equation}
We have
\[
\frac{m_1 m_2 \cdots m_{k+1}}{n_{k+1}} = \frac{m_1 m_2 \cdots m_{k+1}}{n_{k}(m_{k+1}+s_{k+1})}
= \frac{m_1 m_2 \cdots m_k}{n_k} \cdot \frac{m_{k+1}}{m_{k+1}+s_{k+1}}
\]
by construction.  

We may then choose $m_{k+1}$ and $s_{k+1} \neq 0$ to satisfy
\begin{equation}\label{dimbound}
\frac{c}{2} < \frac{m_1 m_2 \cdots m_{k+1}}{n_{k+1}} < \frac{c}{2} + \frac{1}{2^{k+1}},
\end{equation}
whence (\ref{closetodrr}) holds for all $k \in \mathbb{N}$.
Theorem 1 of \cite{V1} shows that the points $x_{i-1}^1,\ldots,x_{i-1}^{s_{i-1}} \in X_{i-1}$,
$i \in \mathbb{N}$, may be chosen in a manner which makes the (unital) limit algebra $A_c = \lim_{i \to \infty}
(B_i,\phi_i)$ simple.  By (\ref{dimbound}) we have
\[
\lim_{i \to \infty} \frac{\mathrm{dim}(X_i)}{n_i} = \lim_{i \to \infty} \frac{m_1 m_2 \cdots m_i}{n_i}
= 2 \left(\frac{c}{2}\right) = c,
\]
whence $\mathrm{drr}(A_c) \leq c$.

In order to conclude that $\mathrm{drr}(A_c) = c$, we must prove that any other decomposition
of $A$ satisfies
\[
\liminf_{i \to \infty} \ \mathrm{max}_{1 \leq l \leq m_i}\left\{\frac{\mathrm{dim}(X_{i,l})}{\mathrm{rank(p_{i,l})}}\right\} \geq c.
\]
To this end we will employ the ordered $\mathrm{K}_0$-group of $A_c$.  

Let $\xi$ denote the Hopf line bundle over $\mathrm{S}^2$, and $\theta_l$ the trivial
vector bundle of complex fibre dimension $l \in \mathbb{N}$ over an arbitrary compact
Hausdorff space $X$.  In \cite{V1} it is proved that the $\mathrm{K}^0(X_i)$-class
\[
y_i:=[\xi^{\times m_1 m_2 \cdots m_i}] - [\theta_1]
\]
is not positive in either $\mathrm{K}^0(X_i)$ or $\mathrm{K}_0(A_c)$.  
By Lemma \ref{onestate}, $\mathrm{K}_0(A_c)$ admits a unique state $s$, which is
realised on $B_i$ as the normalised geometric state --- the state which returns the virtual dimension
of a $\mathrm{K}^0(X_i)$-class divided by $n_i$.  Thus, 
\[
s(y_i) = \frac{m_1 m_2 \cdots m_i-1}{n_i} \stackrel{i \to \infty}{\longrightarrow} \frac{c}{2}.
\]

Suppose that there exists an injective decomposition $A_c = \lim_{i \to \infty}(A_i,\phi_i)$ with 
\[
A_i := \bigoplus_{l=1}^{m_i} p_{i,l} (\mathrm{C}(Y_{i,l}) \otimes \mathcal{K}) p_{i,l}
\]
and such that
\[
\liminf_{i \to \infty} \mathrm{max}_{1 \leq l \leq m_i}\left\{\frac{\mathrm{dim}(Y_{i,l})}{\mathrm{rank(p_{i,l})}}\right\} < c.
\]
By compressing the inductive sequence in this decomposition we may assume that
\[
\lim_{i \to \infty} \mathrm{max}_{1 \leq l \leq m_i}\left\{\frac{\mathrm{dim}(Y_{i,l})}{\mathrm{rank(p_{i,l})}}\right\} < c.
\]
Choose $i_0 \in \mathbb{N}$ and $\epsilon > 0$ such that
\[ 
\mathrm{max}_{1 \leq l \leq m_i}\left\{\frac{\mathrm{dim}(Y_{i,l})}{\mathrm{rank(p_{i,l})}}\right\} < c - \epsilon, \ \forall i \geq i_0.
\]
Choose $j \geq i_0$ large enough so that $s(y_j) > (c - \epsilon)/2$ and $y_j \in \mathrm{K}_0(A_j)$.
Put
\[
y_j^l = y_j|_{\mathrm{K}_0(p_{j,l} (\mathrm{C}(Y_{i,l}) \otimes \mathcal{K}) p_{j,l})}, \ 1 \leq l \leq m_j, 
\]
so that $y_j = \oplus_l y_j^l$.  Applying Lemma \ref{evenstate}, we may have that
\[
s(y_j^l) > \frac{c - \epsilon}{2}, \ 1 \leq l \leq m_j. 
\]
This, in turn, implies that the virtual dimension of each $y_j^l$ is greater than $\lceil \mathrm{dim}(Y_{j,l})/2 \rceil$,
whence each $y_j^l$ is positive in $\mathrm{K}_0(p_{j,l} (\mathrm{C}(Y_{i,l}) \otimes \mathcal{K}) p_{j,l})$.
But then $y_j$ must be positive, contradicting our choice of $y_j$.

That $A_c$ has stable rank one follows from Lemma 9 and Proposition 10 of \cite{V1}, upon
noticing that the general construction of $A_c$ is of the type described in section 2 of
the same paper.

\end{proof} 

\begin{cors}\label{perfrad}
Let $c \in \mathbb{R}^+$.  Then, with $A_c$ as in Theorem \ref{drrrange}, we have
\[
\mathrm{inf} \{s \in \mathbb{R}| \ \mathrm{K}_0 A_c \ \mathrm{has} \ s\mathrm{-strict} \ \mathrm{comparison} \ \} 
= \frac{\mathrm{drr}(A)}{2} = \frac{c}{2}.
\]
\end{cors}

\begin{proof}
$A_c$ has a $\mathrm{K}_0$-group which admits a unique state by Lemma \ref{onestate}.  
We may thus apply Theorem \ref{drrstab} to conclude that $\mathrm{K}_0A_c$ has 
$(\mathrm{drr}(A)/2)$-strict comparison.  This proves the corollary if $c=0$.  

If $c>0$, then $\mathrm{K}_0A_c$
does not have $s$-strict comparison for any $s < c/2$.  
Indeed, the element $y_i = [\xi^{\times m_1 m_2 \cdots m_i}] - [\theta_1]$ is not
positive in $\mathrm{K}_0A_i$, and neither is its image in $\mathrm{K}_0A_c$.
Applying the geometric state on $\mathrm{K}_0A_i$, one has
\[
s([\theta_1]) = \frac{1}{n_i}; \ \ s([\xi^{\times m_1 m_2 \cdots m_i}]) = 
\frac{m_1 m_2 \cdots m_i}{n_i}.
\]
Choosing $i$ large enough so that $c/2 - 1/n_i > s$ we have
\[
s([\theta_1]) + s < s([\xi^{\times m_1 m_2 \cdots m_i}]),
\]
yet $[\theta_1] \nleq [\xi^{\times m_1 m_2 \cdots m_i}]$.  The corollary follows.
\end{proof}

\begin{cors}\label{uncount}
The class of simple, unital and non-$\mathcal{Z}$-stable AH algebras is uncountable.
\end{cors}

\begin{proof}
The algebra $A_c$ of Theorem \ref{drrrange} has a perforated ordered
$\mathrm{K}_0$-group for each $c \neq 0$.  Theorem 1 of \cite{GJS} states that
a simple, unital, finite, and $\mathcal{Z}$-stable $C^*$-algebra has a weakly unperforated
ordered $\mathrm{K}_0$-group, whence the $A_c$s in question are non-$\mathcal{Z}$-stable.
\end{proof}

The pairwise non-isomorphic algebras constructed
in the proof of Theorem \ref{drrrange} are difficult to distinguish from one another
without using the dimension-rank ratio.  Straightforward calculation
shows that, for each $c \neq 0$, $\mathrm{T}(A_c)$ is a
Bauer simplex with extreme boundary homeomorphic to $(S^2)^{\infty}$,
and $\mathrm{K}_1A_c = 0$.  Computing the ordered group $\mathrm{K}_0 A_c$
is not feasible --- the order structure on $\mathrm{K}_0 (S^2)^n$ is not known
for general $n$.

\section{Abstracting the dimension-rank ratio}

The dimension-rank ratio functions well as an invariant tailored
for the study of unital and non-$\mathcal{Z}$-stable AH algebras, so it is
natural to ask whether there exists an invariant defined for any 
unital and stably finite $C^*$-algebra which recovers (or is at least 
closely related to) the dimension-rank ratio upon restricting
to the subclass of unital AH algebras.  In this section, we present a 
candidate for such an invariant.

One could, in light of Corollary \ref{perfrad}, be forgiven for wondering briefly if 
the extended real
\[
\mathrm{inf} \{s| \ \mathrm{K}_0 A \ \mathrm{has} \ s\mathrm{-strict} \ \mathrm{comparison} \ \}
\]
might be the invariant we seek.  The algebra $\mathrm{C}([0,1]^n)$ dispels
this notion:  its $\mathrm{K}_0$-group has comparison, yet $\mathrm{drr}(\mathrm{C}([0,1]^n)) = n$.  
There is, however, a different version of ordered $\mathrm{K}$-theory,
whose prospects for recovering the dimension-rank ratio are distinctly better than those
of the $\mathrm{K}_0$-group.

Let $A$ be a $C^*$-algebra.  We recall the definition of the Cuntz semigroup $W(A)$ from \cite{C}.  (Our 
synopsis is essentially that of \cite{R3}.)  Let
$\mathrm{M}_n(A)^+$ denote the positive elements of $\mathrm{M}_n(A)$,
and let $\mathrm{M}_{\infty}(A)^+$ be the disjoint union $\cup_{i=n}^{\infty} 
\mathrm{M}_n(A)^+$.  For $a \in \mathrm{M}_n(A)^+$ and $b \in \mathrm{M}_m(A)^+$
set $a \oplus b = \mathrm{diag}(a,b) \in \mathrm{M}_{n+m}(A)^+$, and write 
$a \precsim b$ if there is a sequence $\{ x_k \}$ in $\mathrm{M}_{m,n}(A)$ such
that $x_k^* b x_k \rightarrow a$.  Write $a \sim b$ if $a \precsim b$ and 
$b \precsim a$.  Put $W(A) = \mathrm{M}_{\infty}(A)^+ / \sim$, and let
$\langle a \rangle$ be the equivalence class containing $a$.  Then,
$W(A)$ is a positive ordered Abelian semigroup when equipped with the
relations:
\begin{displaymath}
\langle a \rangle + \langle b \rangle = \langle a \oplus b \rangle, \ \ \ \ \
\langle a \rangle \leq \langle b \rangle \Longleftrightarrow a \precsim b, \ \ \ \ \
a,b \in \mathrm{M}_{\infty}(A)^+.
\end{displaymath}
The relation $\precsim$ reduces to Murray-von Neumann comparison when $a$ and $b$
are projections and $A$ is stably finite.

In the case of a stably finite $C^*$-algebra $A$, the Cuntz semigroup 
may be thought of as a generalised version of the semigroup of
Murray-von Neumann equivalence classes of projections in $\mathrm{M}_{\infty}(A)$.
If $A$ is unital, then we scale $W(A)$ with $\langle 1_A \rangle$.  Let $S(W(A))$
denote the set of additive and order preserving maps from $W(A)$ to $\mathbb{R}^+$
having the property that $s(\langle 1_A \rangle) = 1$, $\forall s \in S(W(A))$.
Such maps are called {\it states}.  Given $\tau \in \mathrm{QT}(A)$, one may 
define a map $s_{\tau}:\mathrm{M}_{\infty}(A)^+ \to \mathbb{R}^+$ by
\begin{equation}\label{ldf}
s_{\tau}(a) = \lim_{n \to \infty} \tau(a^{1/n}).
\end{equation}
This map is lower semicontinous, and defines a state on $W(A)$.  Such maps are
called {\it lower semicontinuous dimension functions}, and the set of them 
is denoted $\mathrm{LDF}(A)$. $\mathrm{QT}(A)$ is a simplex (Theorem II.4.4, \cite{BH}), 
and the map from $\mathrm{QT}(A)$ to $\mathrm{LDF}(A)$ defined by (\ref{ldf}) is 
bijective and affine (Theorem II.2.2, \cite{BH}).

\begin{dfs}\label{rc} Let $A$ be an unital and stably finite $C^*$-algebra, and let $r>0$.  
\begin{enumerate}
\item[(i)] Say that $A$ has $r$-comparison if whenever one has positive elements
$a,b \in \mathrm{M}_{\infty}(A)$ such that
\[
s(\langle a \rangle) + r < s(\langle b \rangle), \ \forall s \in \mathrm{LDF}(A),
\]
then $\langle a \rangle \leq \langle b \rangle$ in $W(A)$.
\item[(ii)] Define the radius of comparison
of $A$, denoted $\mathrm{rc}(A)$, to be 
\[
\mathrm{inf} \{r \in \mathbb{R}^+ | \ (W(A), \langle 1_A \rangle) \ \mathrm{has} \ r-\mathrm{comparison} \ \}
\]
if it exists, and $\infty$ otherwise.
\end{enumerate}
\end{dfs}

We summarise some properties of the radius of comparison which are more or less immediate from its
definition.  (Compare with Propositions \ref{drrprop} and \ref{rhomprop}.)

\begin{props}\label{rcprop}
Let $A,B$ be unital and stably finite $C^*$-algebras.  Then:
\begin{enumerate}
\item[(i)] $\mathrm{rc}(A \oplus B) = \mathrm{max} \{\mathrm{rc}(A), \mathrm{rc}(B)\}$;
\item[(ii)] if $k \in \mathbb{N}$, then $\mathrm{rc}(\mathrm{M}_k(A)) = (1/k)\mathrm{rc}(A)$;
\item[(iii)] if $I$ is an ideal of $A$, $\pi:A \to A/I$ is the
quotient map, and
\[
\pi^{\sharp}:\mathrm{QT}(A/I) \to \mathrm{QT}(A)
\]
is surjective, then $\mathrm{rc}(A/I) \leq \mathrm{rc}(A)$;
\item[(iv)] $A$ has $\mathrm{rc}(A)$-(FCQ);
\item[(v)] $(\mathrm{K}_0A^+,[1_A])$ has $\mathrm{rc}(A)$-strict comparison.
\end{enumerate}
\end{props}

\begin{proof}
For (i), use the fact that $\mathrm{LDF}(A \oplus B)$ is the convex hull of $\mathrm{LDF}(A)$ and $\mathrm{LDF}(B)$
to obtain $\mathrm{rc}(A \oplus B) \leq \mathrm{max} \{\mathrm{rc}(A), \mathrm{rc}(B)\}$.  The reverse inequality
will follow from (iii).

(ii) is straightforward from Definition \ref{rc}.

For (iii), let there be given positive elements $a,b \in \mathrm{M}_{\infty}(A/I)$ such that
\[
s(\langle a \rangle) + r < s(\langle b \rangle), \ \forall s \in \mathrm{LDF}(A/I), \ \mathrm{some} \ r > \mathrm{rc}(A).
\]
We may find
positive elements $\tilde{a}, \tilde{b} \in \mathrm{M}_{\infty}(A)$ such that
$\pi(\tilde{a}) = a$ and $\pi(\tilde{b}) = b$ (lift to self-adjoint elements and apply
the functional calculus).  Let $d_{\tau} \in \mathrm{LDF}(A)$ be a state corresponding
to a normalised quasi-trace $\tau$ on $A$.  Then, $\tau = \pi^{\sharp}(\eta)$ for some
$\eta \in \mathrm{QT}(A/I)$ by assumption, and
\[
d_{\tau}(\langle \tilde{a} \rangle) = \lim_{n \to \infty} (\tau(\tilde{a}^{1/n}) =
\lim_{n \to \infty} \eta(a^{1/n}) = s_{\eta}(\langle a \rangle)
\]
for the state $s_{\eta}$ corresponding to some $\eta \in \mathrm{QT}(A/I)$.
It follows that 
\[
d_{\tau}(\langle \tilde{a} \rangle) + r < d_{\tau}(\langle \tilde{b} \rangle), \ \forall d \in \mathrm{LDF}(A), \ \mathrm{some} \ r > \mathrm{rc}(A),
\]
whence $\tilde{a} \precsim \tilde{b}$ in $W(A)$.  This implies the existence of
a sequence $(v_k) \subseteq \mathrm{M}_{\infty}(A)$ such that
\[
v_k^* \tilde{b} v_k \stackrel{k \to \infty}{\longrightarrow} \tilde{a}.
\]
Applying $\pi$ to the expression above shows that $a \precsim b$ in $W(A/I)$, as desired.

(iv) and (v) follow from the fact that there is an order unit preserving order embedding of the scaled ordered
semigroup $(V(A),[1_A])$ of Murray-von Neumann equivalence classes of projections in $\mathrm{M}_{\infty}(A)$
into $(W(A),\langle 1_A \rangle)$ whenever $A$ is stably finite (cf. \cite{R4}).
\end{proof}

The next proposition is the analogue of Proposition \ref{closed} for the radius of
comparison.

\begin{props}
Let $A$ be an unital and stably finite $C^*$-algebra for which every $\tau \in \mathrm{QT}(A)$
is faithful.  Then, the set 
\[
B := \{r \in \mathbb{R}^+ \ | \ W(A) \ \mathrm{has} \ r-\mathrm{comparison} \ \}
\]
is closed.  In other words, $W(A)$ has $\mathrm{rc}(A)$-comparison.
\end{props}

\begin{proof}
If $B = \emptyset$, then it is closed;  suppose that $B \neq \emptyset$.
As in the proof of Proposition \ref{closed}, we need only prove that $\alpha :=
\mathrm{inf}(B) \in B$.  Let there be given $a,b \in W(A)$ satisfying
\[
s(a) + \alpha < s(b), \ \forall s \in \mathrm{LDF}(A).
\]

Suppose first that $a = \langle p \rangle$ for some projection $p \in \mathrm{M}_{\infty}(A)$.
Then, the map $\gamma_a:\mathrm{QT}(A) \to \mathbb{R}^+$ given by $\gamma_a(s) = s(a)$
is continuous.  By \cite[Proposition 2.7]{PT}, the map $\gamma_b:\mathrm{QT}(A) \to \mathbb{R}^+$ 
given by $\gamma_b(s) = s(a)$ is lower semicontinuous.  It follows that $\gamma_b - \gamma_a$
is lower semicontinuous and strictly positive on $\mathrm{QT}(A)$.  Since
$\mathrm{QT}(A)$ is compact, $\gamma_b - \gamma_a$ achieves a lower bound $\delta > 0$,
whence
\[
s(a) + \alpha + \delta/2 < s(b), \ \forall s \in \mathrm{LDF}(A).
\]
$W(A)$ has $(\alpha + \delta/2)$-comparison, and so $a \leq b$, as desired.

Now suppose that $a$ is not Cuntz equivalent to any projection.  By the functional
calculus, we conclude that $0$ is not an isolated point of the spectrum of $a$.
Viewing $a$ as the function $f(t)=t$ on its spectrum, we denote by $(a-\epsilon)_+$
the function $\mathrm{max} \{0,f(t)-\epsilon\}$ on the spectrum of $a$.
By \cite[Proposition 2.6]{KR}, proving that $a \leq b$ is equivalent to 
proving that $\langle (a-\epsilon)_+ \rangle \leq b$, $\forall \epsilon > 0$.
Let $g_\epsilon(t) \in C^*(a)$ be a function supported on $(0,\epsilon) \cap \sigma(a) (\neq \emptyset)$,
where $\sigma(a)$ denotes the spectrum of $a$.  Since $g_\epsilon(t) + (a-\epsilon)_+ \leq f(t) = a$,
we have
\[
s(a) - s((a-\epsilon)_+) \geq s(g_\epsilon), \ \forall s \in \mathrm{LDF}(A).
\]
Let $\mathrm{supp}(\bullet)$ denote the support of a function.
Each $s \in \mathrm{LDF}(A)$ is implemented on $C^*(a)$ by a probability measure
$\mu_s$ in the following sense:  for any $d \in C^*(a)$, $s(d) = \mu_s(\mathrm{supp}(d))$.
Moreover, our assumption about the faithfulness of quasitraces on $A$ implies
that $\mu_s(U) > 0$ for every open subset $U$ of $\sigma(a)$.  Thus, the map
$\gamma_{g_\epsilon}:\mathrm{QT}(A) \to \mathbb{R}^+$ given by $\gamma_{g_\epsilon}(s)
= s(g_\epsilon) = \mu_s( (0,\epsilon) \cap \sigma(a))$ is strictly positive, and,
as above, lower semicontinuous.  It follows that $\gamma_{g_\epsilon}$ achieves a
lower bound on $\mathrm{QT}(A)$, say $\delta_\epsilon$.  Now
\[
s((a-\epsilon)_+) + \alpha + \delta_\epsilon/2 < s(b), \ \forall s \in \mathrm{LDF}(A).
\]
$W(A)$ has $(\alpha + \delta_\epsilon/2)$-comparison for every $\epsilon>0$, whence
$\langle (a-\epsilon)_+ \rangle \leq b$, $\forall \epsilon>0$.
\end{proof}

Recall that $W(A)$ is said to be {\it almost unperforated} if $x \leq y$ in $W(A)$ whenever $mx \leq ny$
for natural numbers $m > n$ (\cite{R4}).  
\begin{props}\label{rczerotoaup}
Let $A$ be an unital and stably finite $C^*$-algebra for which every $\tau \in \mathrm{QT}(A)$
is faithful.  If $\mathrm{rc}(A)=0$, then $W(A)$ is almost unperforated.
\end{props}

\begin{proof}
Let $m>n$ be natural numbers, and $x,y \in W(A)$ such that $mx \leq ny$.  
For any $s \in \mathrm{LDF}(A)$ we have the following string of inequalities:
\begin{eqnarray*}
0 & \leq & n \cdot s(y) - m \cdot s(x) \\
0 & \leq & n(s(y)-s(x)) - (m-n)s(x) \\
\frac{(m-n)s(x)}{n} & \leq & s(y)-s(x).
\end{eqnarray*}
The map $\gamma:\mathrm{QT}(A) \to
\mathrm{R}^+$ given by $s \mapsto s(x)$ is thus strictly positive (since each
$\tau \in \mathrm{QT}(A)$ is faithful) and lower
semicontinuous (\cite[Proposition 2.7]{PT}).  
Since $\mathrm{QT}(A)$ is compact, $\gamma$ achieves a minimum
value $\delta > 0$.  Now $s(x) + \delta/2 < s(y)$, $\forall s \in \mathrm{LDF}(A)$.  Since 
$\mathrm{rc}(A) = 0$, $W(A)$ has $(\delta/2)$-comparison.  We conclude that 
$x \leq y$ in $W(A)$, as desired.
\end{proof}

\begin{thms}[R{\o}rdam, Corollary 4.6, \cite{R4}]\label{auptorczero}
Let $A$ be a simple, unital, exact, stably finite $C^*$-algebra with $W(A)$ almost unperforated.
Then, $W(A)$ has $0$-comparison, and $\mathrm{rc}(A)=0$.
\end{thms}

\noindent
Combining Proposition \ref{rczerotoaup}, Theorem \ref{auptorczero}, and Corollary
\ref{taf}, we conclude that $\mathrm{drr}=0$ and $\mathrm{rc}=0$ are equivalent
for simple and infinite-dimensional AH algebras of real rank zero and stable rank one.
We shall see in Corollary \ref{dimbd} below that if a semi-homogeneous 
algebra $A$ has $W(A)$ almost unperforated and  
spectrum a CW-complex, then the dimension of its spectrum,
and hence its dimension-rank ratio, is at most four;
by Proposition \ref{rczerotoaup}, this conclusion holds {\it a fortiori} 
if the said semi-homogeneous algebra has $\mathrm{rc}=0$.  If every 
finite-dimensional representation of $A$ is large, then 
$\mathrm{rc}=0$ implies that $\mathrm{drr} \approx 0$.
$\mathrm{drr}(A)=0$ implies that the spectrum of $A$ is zero-dimensional;  $W(A)$ is
then almost unperforated by Theorem 3.4 of \cite{P}.
Taken together, these results show the condition $\mathrm{rc}=0$ to be an appropriate 
abstraction of the condition $\mathrm{drr}=0$.

\begin{thms}\label{rcbound}
Let $X$ be a CW-complex of finite dimension $n$, $p \in \mathrm{C}(X) \otimes \mathcal{K}$ a projection, 
and $m$ the greatest nonnegative integer such that $2m<n$.  Then,
\[
\mathrm{rc}(p(\mathrm{C}(X) \otimes \mathcal{K})p) \geq \frac{m-1}{\mathrm{rank}(p)}.
\]
\end{thms}

\begin{proof}
The theorem is trivial if $m \leq 1$, so suppose that $m \geq 2$.
Choose an $n$-cell of $X$, say $E$.  There is a subset $A$ of $E^{\circ}$
homeomorphic to $(-1,1)^n$.  Let $\psi:A \to (-1,1)^{2m+1}$ be the projection
onto the first $2m+1$ co-ordinates of $A$, and let $d$ be the usual Euclidean
metric on $\mathrm{Im}(\psi) = (-1,1)^{2m+1}$.  Put
\[
Y:=\{ (x_1,\ldots,x_{2m+1}) \in \mathrm{Im}(\psi) \ | \ d \left( (x_1,\ldots,x_{2m+1}),(0,\ldots,0) \right) = 1/2 \}  
\]
and
\[
S:=\{ (x_1,\ldots,x_{2m+1}) \in \mathrm{Im}(\psi) \ | \ 1/3 < d \left( (x_1,\ldots,x_{2m+1}),(0,\ldots,0) \right) < 2/3 \}.
\]
Let $r:S \to Y$ be the projection along rays emanating from $(0,\ldots,0) \in \mathrm{Im}(\psi)$.
Put $O = \psi^{-1}(S)$ and $\pi = r \circ \psi$.  We now have a closed 
subset $Y$ of $E^{\circ}$ homeomorphic to $S^{2m}$, an open set $O$ such that
$E^{\circ} \supseteq O \supseteq Y$, and a continuous map $\pi:O \to Y$
such that $\pi|Y = \mathrm{id}_Y$.

Recalling the description of $(\mathrm{K}^0S^{2m}, \mathrm{K}^0{S^{2m}}^+)$ 
from the example following Proposition \ref{rhomprop},
let $\xi_m$ be a complex vector bundle over $Y$ whose $\mathrm{K}^0$-class corresponds
to $m \oplus 1 \in \mathbb{Z} \oplus \mathbb{Z} \cong \mathrm{K}^0S^{2m}$.
$\xi_m$ can be realised inside $M_{2m}(\mathrm{C}(X))$.
If $\theta_1$ is the trivial complex line bundle over $Y$, then
the class $[\theta_1]$ corresponds to the element $1 \oplus 0 \in \mathrm{K}^0S^{2m}$ and is
clearly not dominated by $[\xi_m]$.  Let $f:X \to [0,1]$ be a continuous
function which vanishes off $O$ and takes the value $1$ at every point in the closure of some
open set $V \supseteq Y$ such that $\overline{V} \subseteq O$.
Define positive functions $a, b \in \mathrm{M}_{2m}(\mathrm{C}(X))$ by
\[
a(x) = f(x)  \pi^*(\xi_m); \ \ b(x) = f(x) \pi^*(\theta_1).
\]
We may think of $a$ and $b$ as being contained in $\mathrm{M}_k(p(\mathrm{C}(X) \otimes
\mathcal{K})p)$ for some sufficiently large $k \in \mathbb{N}$.  

We claim that $\langle b \rangle \nleq \langle a \rangle$ in $W(p(\mathrm{C}(X) \otimes
\mathcal{K})p)$.  Indeed, since $f(y) =1$, $\forall y \in Y$, and $\pi|Y = \mathrm{id}_Y$,
we have that
\[
a(y) = \xi_m(y), \ b(y) = \theta_1(y), \ \forall y \in Y.
\]
$\langle b \rangle \leq \langle a \rangle$ implies that
$\langle b|Y \rangle \leq \langle a|Y \rangle$ in $W(\mathrm{C}(Y))$, but 
the second inequality contradicts the fact that $\theta_1$ is not Murray-von
Neumann equivalent to a subprojection of $\xi_m$ (remember that the Cuntz equivalence
relation reduces to Murray-von Neumann equivalence on projections in a stably finite
algebra).  The claim follows.

Choose a continuous function $g:X \to [0,1]$ such that $g$ is identically zero on $Y$,
and identically one on the complement of $V$.  Define a positive element $v := g \cdot 
\theta_{n}$.  Since $v$ is zero on $Y$, the argument of the preceding paragraph shows
that 
\[
\langle b \rangle \nleq \langle a \oplus v \rangle.
\]

The lower semicontinuous dimension functions on $A = p(\mathrm{C}(X) \otimes \mathcal{K})p$
correspond to normalised traces on $A$.  This correspondence may be viewed 
as follows:  each normalised trace $\tau$ corresponds to a probability measure
$\mu_{\tau}$ on $X$, and the dimension function $d_{\tau}$ is given by
\[
d_{\tau}(\langle a \rangle) = \int_X \frac{\mathrm{rank}(a)(x)}{\mathrm{rank}(p)} d\mu_{\tau}.
\]

Let $\tau \in \mathrm{T}A$ be given. We have
\begin{eqnarray*}
\mathrm{rank}(p) \cdot d_{\tau}(\langle a \oplus v \rangle) & = & \int_X \mathrm{rank}(a \oplus v)(x) d\mu_{\tau} \\
& = & \int_{X \backslash V} n + \mathrm{rank}(a)(x) d\mu_{\tau} + \int_{V \backslash Y} \mathrm{rank}(v)(x)
+ m d\mu_{\tau} \\
& & \hspace{5mm} + \int_Y m d\mu_{\tau} \\
& \geq & n\mu_{\tau}(X \backslash V) + m\mu_{\tau}(V \backslash Y) + m\mu_{\tau}(Y) \\
& \geq & m
\end{eqnarray*}
and
\[
\mathrm{rank}(p) \cdot d_{\tau}(\langle b \rangle)  =  \int_X \mathrm{rank}(b)(x) d\mu_{\tau} 
= \int_{O} d\mu_{\tau} 
\leq 1.
\]
Thus, for any $s \in \mathrm{LDF}(A)$ we have
\[
s(\langle b \rangle) + \frac{m-1}{\mathrm{rank}(p)} \leq s(\langle a \oplus v \rangle)
\]
while $\langle b \rangle \nleq \langle a \oplus v \rangle$.  The proposition follows.
\end{proof}

The lower bound on $\mathrm{rc}(p(\mathrm{C}(X) \otimes \mathcal{K})p)$ in Theorem \ref{rcbound} is 
close to $\mathrm{drr}(A)/2$, particularly when $\mathrm{dim}(X)$ and $\mathrm{rank}(p)$ are
large.  If a simple unital AH algebra $B$ has $\mathrm{drr}(B) > 0$, then the dimensions
of the spectra of its building blocks and the ranks of the units of these building 
blocks must tend toward infinity, regardless of the injective decomposition chosen. Thus, the bound
of Theorem \ref{rcbound} applied to these building blocks will be all but equal to 
one half of their respective dimension-rank ratios.  One can obtain a lower bound
in the spirit of Theorem \ref{rcbound} for the algebras of Theorem \ref{drrrange}.

The proof of Theorem \ref{rcbound} yields:
\begin{cors}\label{dimbd}
Let $A$ be a semi-homogeneous $C^*$-algebra with spectrum a CW-complex.  If 
$W(A)$ is almost unperforated, then the dimension of the spectrum of $A$ is at most four.
\end{cors}

\begin{proof}  We prove the contrapositive.  Retain the notation used in the proof of Theorem \ref{rcbound}.
Suppose that the dimension of the spectrum of $A$ is at least five.  Construct
$a$ and $b$ as in the proof of Theorem \ref{rcbound}, and notice that 
$a = \pi^*(\xi_2)$.  Theorem \ref{kstab} shows that 
\[
\xi_2 \oplus \xi_2 \oplus \xi_2 \cong \theta_4 \oplus \eta
\]
for some complex vector bundle $\eta$ over $Y$.  In other words, there is a partial 
isometry $v \in \mathrm{M}_{\infty}(\mathrm{C}(Y))$ such that
\[
v^*(\xi_2 \oplus \xi_2 \oplus \xi_2)v = \theta_4.
\]
Let $(g_k)$ be a self-adjoint approximate unit for $\mathrm{C}(O)$.  Put $w_k = g_k \cdot \pi^*(v)$.
Then,
\begin{eqnarray*}
w_k^*(a \oplus a \oplus a)w_k & = & g_k \pi^*(v^*(\xi_2 \oplus \xi_2 \oplus \xi_2)v) g_k \\
& = & g_k \pi^*(\theta_4) g_k \\
& = & g_k (\oplus_{j=1}^4 \pi^*(\theta_1)) g_k \\
& = & \oplus_{j=1}^4 g_k b g_k \\
& \stackrel{k \to \infty}{\longrightarrow} & \oplus_{j=1}^4 b.
\end{eqnarray*} 
This is precisely the statement that $4 \langle b \rangle \leq 3 \langle a \rangle$.  
$\langle b \rangle \nleq \langle a \rangle$ by the proof of Theorem \ref{rcbound}, 
and the corollary follows.
\end{proof}  

\begin{props}\label{rcbound2}
For any $r \in \mathbb{R}^+ $, there is a simple unital AH algebra $A$ such
that $\mathrm{rc}(A) \geq r = \mathrm{drr}(A)/2$.  
\end{props}

\begin{proof}
There is nothing to prove when $r=0$, so fix $r >0$.  For a $C^*$-algebra $A$, let $V(A)$ denote
the semigroup of  Murray-von Neumann equivalence classes
of projections in $\mathrm{M}_{\infty}(A)$.  The algebra $A_{2r}$
constructed in the proof of Theorem \ref{drrrange} has $\mathrm{drr} = 2r$ and stable rank one.
It follows that there is an order unit preserving order isomorphism
\[
(\mathrm{K}_0(A_{2r})^+,[1_{A_{2r}}]) \cong (V(A_{2r}),[1_{A_{2r}}]).
\]
Since $A_{2r}$ is stably finite,
there is an order unit preserving order embedding of $(V(A_{2r}),[1_{A_{2r}}])$
into $(W(A_{2r}),\langle 1_{A_{2r}} \rangle)$.  The proof of Corollary \ref{perfrad}
shows that for any $t < r$, there are projections $p,q \in \mathrm{M}_{\infty}(A)$
such that $[p] \nleq [q]$, yet $s(p) + t < s(q)$ for the unique state $s \in S(\mathrm{K}_0(A_{2r}),
[1_{A_{2r}}])$.  

Let $d_{\tau} \in \mathrm{LDF}(A_{2r})$ be induced by $\tau \in \mathrm{T}(A)$.
$d_{\tau}$ gives rise to an element of $S(\mathrm{K}_0(A_{2r})$, and so agrees with $s$
on the image of $\mathrm{K}_0(A_{2r})^+$ in $W(A_{2r})$.  In particular,
\[
d_{\tau}(\langle p \rangle) +t < d_{\tau}(\langle q \rangle), \ \forall d_{\tau} \in \mathrm{LDF}(A_{2r}).
\]
The existence of an order unit preserving order embedding 
\[
\iota: (\mathrm{K}_0(A_{2r})^+,[1_{A_{2r}}]) \to (W(A_{2r}),\langle 1_{A_{2r}} \rangle)
\]
implies that $\langle p \rangle$ is not less than $\langle q \rangle$ in $W(A_{2r})$, whence $\mathrm{rc}(A_{2r}) \geq t$;
$t$ was arbitrary, and the proposition follows.  
\end{proof}

One wants an upper bound on the radius of comparison of $A = p(\mathrm{C}(X) \otimes \mathcal{K})p$
of the form
\begin{equation}\label{upperbound}
\mathrm{rc}(A) \leq K \mathrm{drr}(A), \ K>0,
\end{equation}
where $X$ is a CW-complex, $p \in \mathrm{C}(X)$ is a projection, and $K$ 
is independent of our choice of $X$ and $p$.  (This bound holds already in
the case $\mathrm{drr}(A)=0$ by Theorem 3.4 of \cite{P}.)  This would complete the confirmation of
the radius of comparison as the correct abstraction of the dimension-rank
ratio.  Applied to the algebras of Theorem \ref{rcbound},
it would show that the radius of comparison roughly determines the dimension rank
ratio. 
Philosophically, asking for the bound in (\ref{upperbound})
is reasonable --- it amounts to asking for stability properties in the Cuntz
semigroup analogous to the the stability properties of vector bundles (cf. Theorem \ref{kstab}):
\begin{qus}\label{wstab}
Does there exist a constant $K>0$ such that for any compact Hausdorff space $X$ and
any positive elements $a,b \in \mathrm{M}_{\infty}(\mathrm{C}(X))$ satisfying
\[
\mathrm{rank}(b)(x) - \mathrm{rank}(a)(x) \geq K \mathrm{dim}(X), \ \forall x \in X,
\] 
one has $a \precsim b$ in $W(\mathrm{C}(X))$?
\end{qus}

\noindent
It follows more or less directly from Theorem \ref{kstab}, (ii), that Question \ref{wstab} has a positive answer
upon restricting to positive elements whose rank functions take at most two values, one of 
which is zero, but this partial result does not address the essential difficulties
of the question.  Nevertheless, an affirmative answer seems likely.  
To generate interest in Question \ref{wstab}, we outline an application of a positive answer
to it.
\begin{conjs}\label{ce}
There exists a simple, unital, separable, and nuclear $C^*$-algebra $A$ of stable rank one such
that 
\[
A \ncong \mathrm{M}_n(A), \ \mathrm{some} \ n \in \mathbb{N},
\]
yet
\[
(V(A),[1_A]) \cong (V(\mathrm{M}_n(A)),[1_{\mathrm{M}_n(A)}]) \cong (\mathbb{Q}^+,1).
\]
\end{conjs}

\noindent
The algebras $A$ and $\mathrm{M}_n(A)$ thus constitute a particularly strong counterexample to 
Elliott's classification conjecture for simple, separable, and nuclear $C^*$-algebras (cf. \cite{R1}).

\vspace{2mm}
\noindent
{\it Sketch of proof.}  The proof of Theorem 1.1, \cite{To2}, contains a construction of a simple
unital AH algebra $A$ of stable rank one which has the properties that $\mathrm{rc}(A) > 1/2$, and
\[
(V(A),[1_A]) \cong (\mathbb{Q}^+,1).
\]  
Explicitly (cf. \cite{To2}), one has 
\[
A = \lim_{i \to \infty}(\mathrm{M}_{n_i}(\mathrm{C}([0,1]^{m_i})),\phi_i),
\]
where $m_i \leq n_i$.  If Question \ref{wstab} has a positive answer, then we may conclude that
$\mathrm{rc}(A) \leq K$.  Choose $n > 2K$.  Then, 
\[
(V(\mathrm{M}_n(A)),[1_{\mathrm{M}_n(A)}]) \cong (\mathbb{Q}^+,1),
\]   
but $\mathrm{rc}(\mathrm{M}_n(A)) < 1/2$ by part (ii) of Proposition \ref{rcprop}.  It follows that
$A \ncong \mathrm{M}_n(A)$, as desired. \hfill $\Box$

\end{document}